\documentclass[a4paper,fleqn]{cas-sc}

\usepackage[authoryear,longnamesfirst]{natbib}
\usepackage{amsmath}
\usepackage{amssymb}
\usepackage{booktabs}
\usepackage{array}
\usepackage{enumitem}
\usepackage{graphicx}
\usepackage{url}
\usepackage{listings}
\usepackage{float}
\usepackage{soul}
\renewcommand{\hl}[1]{#1}
\ExplSyntaxOn
\IfFileExists{thumbnails/cas-email.jpeg}{}
  {\bool_gset_true:N \g_stm_nologo_bool}
\cs_set:Npn \__first_footerline:
{
  \group_begin:
  \small
  \sffamily
  \ifnum\theblind>0\relax
  \else
    \__short_authors: :~
  \fi
  { \rmfamily \itshape Preprint }
  \group_end:
}
\ExplSyntaxOff

\sethlcolor{yellow} 
\def\tsc#1{\csdef{#1}{\textsc{\lowercase{#1}}\xspace}}
\tsc{PSE}
\tsc{MILP}
\tsc{MIQP}
\tsc{PWL}
\tsc{NLP}
\tsc{OMLT}

\begin{document}
\let\WriteBookmarks\relax
\def\floatpagepagefraction{1}
\def\textpagefraction{.001}

\shorttitle{PWLNet for linearized neural surrogates}
\shortauthors{Sanyal et~al.}

\title[mode=title]{PWLNet: A Framework for Linearized Neural Surrogates in Nonlinear Process Optimization}


\author[aff1]{Sunoy Sanyal}[
    orcid=0009-0001-2004-3863
]
\ead{sanyas@rpi.edu}
\credit{Conceptualization, Methodology, Software, Data curation, Formal analysis, Investigation, Visualization, Writing -- original draft, Writing -- review and editing}

\author[aff2, aff3]{Canberk Demir}[
    orcid=0009-0006-5918-4053
]
\ead{cdemir21@ku.edu.tr}
\credit{Conceptualization, Methodology, Software, Data curation, Formal analysis, Investigation, Visualization, Writing -- original draft, Writing -- review and editing}

\author[aff2]{Metin T\"urkay}[
    orcid=0000-0003-4769-6714
]
\ead{mturkay@ku.edu.tr}
\credit{Supervision, Methodology, Resources, Project administration, Writing -- review and editing}

\author[aff4]{Erdal Ayd{\i}n}[
    orcid=0000-0003-3790-5615
]
\ead{erdal.aydin@aalto.fi}
\cormark[1]
\credit{Supervision, Methodology, Resources, Project administration, Writing -- review and editing}
\cortext[cor1]{Corresponding author.}

\affiliation[aff1]{
    organization={Department of Industrial and Systems Engineering, Rensselaer Polytechnic Institute},
    addressline={110 8th St},
    city={Troy, NY},
    postcode={12180},
    country={USA}
}

\affiliation[aff2]{
    organization={Department of Industrial Engineering, Ko\c{c} University},
    addressline={Rumelifeneri Yolu},
    city={Istanbul},
    postcode={34450},
    country={T\"urkiye}
}

\affiliation[aff3]{
    organization={Department of Chemistry, Ko\c{c} University},
    addressline={Rumelifeneri Yolu},
    city={Istanbul},
    postcode={34450},
    country={T\"urkiye}
}

\affiliation[aff4]{
    organization={Department of Chemical and Metallurgical Engineering, Aalto University},
    addressline={Otakaari 1},
    city={Espoo},
    postcode={02150},
    country={Finland}
}

\begin{abstract}
Neural networks can represent nonlinear process systems accurately, but embedding them in mathematical optimization remains difficult because smooth activations introduce nonlinear and nonconvex constraints. PWLNet is a Python framework for constructing optimization-compatible neural surrogates through exact and piecewise-linear reformulations. It combines baseline and physics-informed training, optimization-aware pruning, and Pyomo embeddings based on ReLU mixed-integer linear models, exact tanh nonlinear models, PWL-tanh formulations, and binary-free SOSX branch-and-LP. \hl{PWL-native training retains the same piecewise-linear activation during training and optimization embedding, avoiding additional activation linearization error.} The framework is evaluated using Boston housing dataset, the Williams--Otto process, the Van de Vusse reactor, and a hydroformylation semi-batch reactor. The benchmarks compare predictive accuracy, formulation size, solution time, constraint satisfaction, and downstream decision quality. On Boston housing, SOSX training produced the best common-architecture model, with a test MAE of $1.984$ and $R^2=0.925$. On Williams--Otto, pruning 15\% of the hidden neurons reduced MILP solution time by a factor of 6.5 while preserving the selected operating point. In the Van de Vusse case, a physics-informed ReLU RNN was similarly pruned from 32 to 12 hidden units, enabling tractable exact MILP embedding and a 4\% true-plant economic gain over the best ramp-and-hold policy. For the inverse-response Van de Vusse reactor, terminal-aware shrinking-horizon prediction captured the required build-up-and-conversion dynamics and reached 99.7\% of the direct-shooting reference. \hl{For the hydroformylation reactor, one-step exact-tanh IPOPT reached 99.2\% of the certified reference product with a peak pressure of $20.002$\,bar against the 20-bar limit and $14.36$\,s of accumulated solve time over 280 moves. ReLU and PWL-native controllers achieved comparable product recovery, with ReLU requiring the least solve time. Longer horizons and physics-informed training did not consistently improve closed-loop performance, emphasizing the need to assess realized plant constraints alongside product recovery and computational cost. Overall, these results demonstrate that compact linearized neural surrogates can support accurate and tractable optimization across steady-state and dynamic process applications, while identifying constraint fidelity as an important direction for further improvement.}

\end{abstract}

\begin{graphicalabstract}
\includegraphics[width=\textwidth]{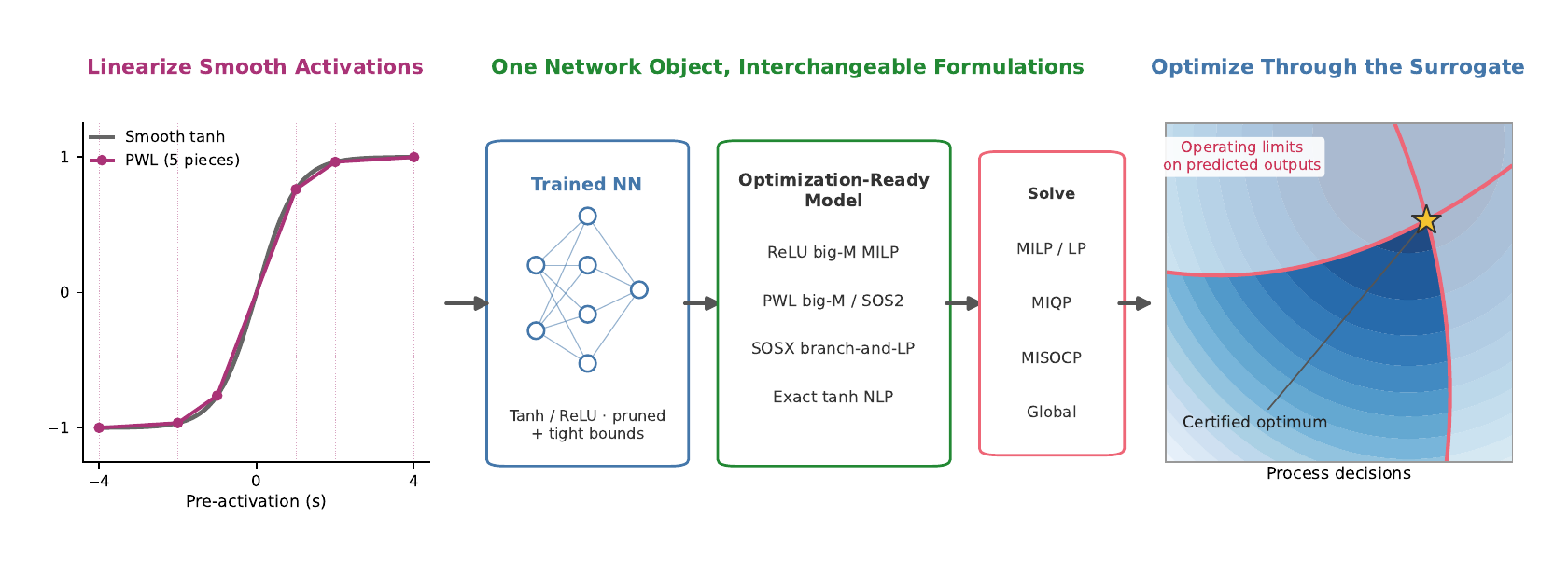}
\end{graphicalabstract}

\begin{highlights}
\item PWLNet provides open-source Python software for optimization-compatible surrogates.
\item ReLU, exact tanh, and piecewise-linear tanh formulations are compared in Pyomo.
\item Pruning and physics-informed training to reduce embedded optimization model size.
\item Improved SOSX algorithm supporting piecewise-linear tanh neural network training.
\item Benchmarks evaluate prediction error, tractability, and downstream decision quality.
\end{highlights}


\maketitle


\section{Introduction}
\label{sec:introduction}

Mathematical optimization and model-based control are central tools in process systems engineering. When accurate mechanistic models are available, classical optimization and control methods provide reliable and interpretable decision-making frameworks. Similarly, when process behavior can be represented adequately using linear or mildly nonlinear approximations, reduced-order models, response surfaces, and conventional data-driven models often provide sufficient accuracy with limited computational burden. In such cases, simpler models remain preferable because they are easier to interpret, easier to validate, and easier to embed in optimization problems.

However, many chemical and biochemical processes exhibit nonlinear, dynamic, and constrained behavior that is difficult to represent using simple model structures. Reaction kinetics, phase equilibrium, transport limitations, saturation effects, recycle streams, and changing operating regimes can produce input--output relationships that are highly nonlinear over the relevant operating region. In addition, first-principles models may be unavailable, incomplete, expensive to simulate, or difficult to maintain when plant conditions change. These limitations create a need for flexible data-driven surrogate models that can approximate complex process behavior while remaining suitable for optimization and control.

Neural networks have therefore received increasing attention as surrogate models for nonlinear process systems. Their ability to approximate complex functions from process data makes them attractive for steady-state optimization, dynamic optimization, model predictive control, fault detection, and soft sensing. More recently, recurrent neural networks and attention-based sequence models have further expanded the ability of data-driven models to represent time-dependent process behavior. These developments suggest that neural networks can serve as powerful process models when mechanistic models are unavailable or insufficiently accurate.

Despite these advantages, the use of neural networks in optimization introduces a fundamental challenge. A trained neural network is straightforward to evaluate as a predictor, but it is not automatically easy to optimize over. In many process applications, the goal is not only to predict the process output for a given input, but to determine the input or control policy that maximizes production, minimizes cost, tracks a target trajectory, or satisfies additional operating constraints. This requires embedding the neural network inside a mathematical optimization model. \hl{For smooth activations such as tanh, direct embedding generally yields a nonconvex nonlinear program, or a mixed-integer nonlinear program when discrete process decisions are also present. In a full-space formulation, intermediate neuron variables and nonlinear activation equations enlarge the optimization model, potentially making global solution computationally demanding as network size increases. This motivates both specialized global optimization methods for smooth neural networks and piecewise-linear representations that enable mixed-integer linear optimization. PWLNet supports the latter approach while retaining exact nonlinear embeddings for comparison.}

A major line of recent work has therefore focused on reformulating trained neural networks into optimization-compatible mathematical programs. ReLU networks are particularly attractive in this context because each activation can be represented using mixed-integer linear constraints. This allows trained ReLU surrogates to be embedded in mixed-integer linear programs and solved using mature MILP solvers \citep{grimstad2019relu_mip,anderson2020strong_relu,tsay2021partition_relu}. Piecewise-linear approximations of smooth activations, such as tanh, provide another route by replacing nonlinear activation functions with linear segments and special ordered set or big-M formulations \citep{koksal2023physics_pwl,koksal2026sosx_training}. These formulations can preserve more of the smooth saturating structure of chemical process responses while maintaining compatibility with mathematical programming solvers.

However, optimization tractability remains strongly dependent on network size. As the process becomes more nonlinear, more dynamic, or higher dimensional, larger neural networks may be needed to achieve sufficient predictive accuracy. Larger networks, in turn, produce larger mathematical programming formulations with more variables, constraints, and, in many cases, binary decisions. This can increase solution time substantially and may render the resulting optimization problem intractable for repeated use in control. Therefore, optimization-aware surrogate modeling requires more than training an accurate neural network. It also requires selecting an appropriate activation function, controlling formulation size, pruning unnecessary network components, tightening variable bounds, and validating the optimized decisions against the original process model or plant simulator.

Despite progress in neural-network reformulation and embedding, existing tools and studies typically address only individual components of the workflow, such as neural-network training, ReLU embedding, algebraic modeling, pruning, or process-specific surrogate optimization. As a result, users often need to manually connect data scaling, network export, formulation construction, pruning, solver configuration, and exact-process validation for each new case study. This makes systematic comparison of activation functions, formulations, solvers, and downstream decision quality difficult and limits reproducibility. The objective of this work is therefore not to propose a single new neural-network formulation or solver algorithm, but to provide a reusable open-source framework that connects existing and newly implemented optimization-compatible surrogate modeling components in a unified workflow.

This work presents PWLNet, an open-source Python framework for optimization-aware neural-network surrogate modeling in process systems engineering. The goal of PWLNet is to connect the full workflow from data generation and neural-network training to model compression, mathematical reformulation, solver selection, and downstream process validation. The framework supports ReLU and tanh neural networks, exact tanh nonlinear formulations, ReLU mixed-integer formulations, piecewise-linear tanh formulations, SOSX-style branch-and-bound workflows, physics-informed training, and optimization-aware pruning. For piecewise-linear tanh formulations, users can specify the number and placement of breakpoints, allowing direct control over the trade-off between approximation accuracy and optimization model size. By exposing these alternatives through a unified interface, PWLNet enables users to compare surrogate accuracy, formulation size, solution time, and decision quality without implementing separate optimization models for each neural-network architecture.

A central motivation of the framework is that no single surrogate formulation is universally superior. ReLU networks may be effective for processes that can be represented compactly using piecewise-affine regions, whereas tanh networks may provide more compact representations for smooth, saturating, or sigmoidal process responses. Similarly, exact nonlinear formulations, MILP formulations, MIQP formulations, conic reformulations, and SOSX-style workflows can exhibit different computational behavior depending on the objective function, constraints, network size, and solver. PWLNet is designed to make these alternatives interchangeable so that formulation choice becomes an empirical and application-specific decision rather than a fixed modeling assumption.

The library also includes tools for reducing optimization model size. Optimization-aware pruning methods, including dead-neuron removal, magnitude-based pruning, neuron pruning, group-lasso-style compression, and Optimal Brain Damage-inspired pruning, are used to remove unnecessary parameters or inactive neurons while preserving surrogate accuracy \citep{lecun1989obd,han2015deep_compression}. Physics-informed training can be used to incorporate known process constraints, conservation relationships, or operating limits during training, improving extrapolation behavior and reducing the network size needed to represent physically meaningful behavior \citep{koksal2023physics_pwl}. These tools are included not only to improve prediction metrics, but to reduce the size and difficulty of the mathematical programs obtained after reformulation.

In addition to conventional trained-network embedding, PWLNet includes a fixed-readout SOSX-MAE training workflow for single-hidden-layer piecewise-linear tanh networks. In this setting, the readout layer is fixed and the remaining training problem is represented using continuous special ordered set variables. SOS2 validity is enforced through branching rather than through binary variables, so each branch-and-bound node can be solved as a linear program. This provides a globally optimized solution for the restricted PWL-tanh training problem and offers a useful benchmark for comparing standard gradient-based training, physics-informed training, and optimization-aware refinement.

The framework is demonstrated on representative process systems engineering case studies, including steady-state process optimization, dynamic nonlinear reactor control, and surrogate-based nonlinear model predictive control. The benchmark suite is chosen to span increasing levels of difficulty: a standard regression problem for baseline predictive performance, a steady-state nonlinear process optimization problem, a dynamic reactor control problem, and a realistic nonlinear process optimization case with more complex dynamics and constraints. These examples are designed to evaluate surrogates not only by test-set prediction error, but also by optimization tractability, formulation size, solver behavior, physical feasibility, exact-process performance, and downstream decision quality. In this way, PWLNet addresses the practical gap between neural networks as predictive models and neural networks as reliable components of mathematical optimization and control workflows.

The main contributions of this work are:
\begin{enumerate}[leftmargin=*]
\item a modular open-source Python framework for optimization-compatible neural surrogate modeling in process systems engineering;
\item an interchangeable workflow for training, embedding, and comparing ReLU, exact tanh, and \hl{PWL-tanh surrogates, including both post-training activation approximation and PWL-native training with a consistent activation representation across prediction and optimization};
\item an integrated implementation of established optimization-aware pruning
and compression methods prior to mathematical model embedding;
\item support for mixed-integer linear, mixed-integer quadratic, conic, nonlinear, and SOSX-style optimization formulations under a unified interface;
\item a benchmark suite covering regression, steady-state process optimization, dynamic reactor behavior, and realistic nonlinear process optimization;
\item exact-process validation of surrogate-optimized decisions to distinguish surrogate objective quality from true process performance;
\item a reproducible evaluation workflow for comparing prediction error, optimization tractability, formulation size, solver performance, and downstream process decision quality.
\end{enumerate}

The remainder of the paper is organized as follows. Section~\ref{sec:background} reviews related work on neural-network surrogates, optimization-compatible reformulations, pruning, physics-informed learning, and learned surrogate control. Section~\ref{sec:framework} describes the overall PWLNet framework. Section~\ref{sec:formulations} introduces the neural-network reformulations considered in this study, and Section~\ref{sec:pruning} describes the optimization-aware pruning workflow. Section~\ref{sec:benchmarks} presents the benchmark case studies and discusses the comparative results. Section~\ref{sec:conclusion} concludes the paper. Section~\ref{sec:software} then describes the repository, licensing, release information, and reproducibility, followed by Data Availability.

\section{Background and Related Work}
\label{sec:background}

\subsection{Surrogate Models in Process Optimization}

Surrogate modeling has a long history in engineering optimization and process systems engineering. Classical surrogate models, including response-surface models, polynomial regression, radial basis functions, Kriging, Gaussian processes, and reduced-order models, have been widely used to approximate expensive simulations and experimental systems \citep{queipo2005surrogate, jones1998ego, rasmussen2006gp}. These methods are especially useful when first-principles models are computationally expensive, unavailable, or inconvenient to embed directly in an optimization problem. In surrogate-based optimization, the surrogate replaces or approximates the original process model, allowing repeated model evaluations during design, scheduling, or control.

Although classical surrogate models remain valuable, they can become difficult to use for high-dimensional, strongly nonlinear, constrained, or dynamic systems. Chemical and biochemical processes often exhibit nonlinear behavior arising from reaction kinetics, phase equilibrium, transport limitations, saturation effects, and operating constraints. These phenomena can produce objective landscapes and state trajectories that are difficult to represent accurately with low-order polynomial or globally smooth response surfaces. Neural networks provide a more flexible surrogate class because they can approximate complex nonlinear mappings from data without requiring the user to prescribe a specific functional form.

Neural-network surrogates have therefore received increasing attention in process systems engineering. They have been applied to steady-state process optimization, dynamic simulation, soft sensing, fault detection, and model predictive control. However, using a neural network as a predictor is fundamentally different from using it as a component of a mathematical optimization problem. Prediction requires only forward evaluation of the trained model, whereas optimization requires searching over the network inputs while satisfying process constraints and operating limits. This distinction motivates the development of optimization-compatible neural-network formulations.

\subsection{Embedding Neural Networks in Mathematical Programs}
\label{sec:embedding}
Embedding trained neural networks into algebraic optimization models has become an important research area. In this setting, the trained network is treated as a surrogate function
\begin{equation}
y = f_{\theta}(x),
\end{equation}
where $x$ denotes process inputs, states, disturbances, or decision variables, $y$ denotes predicted process outputs, and $\theta$ denotes fixed trained network parameters. The optimization problem then searches over $x$ and $y$ subject to the neural-network equations and additional process constraints. The resulting problem class depends strongly on the activation function and objective.

\hl{For smooth neural networks, }\citet{schweidtmann2019embedded} \hl{address the computational burden of full-space embeddings through a reduced-space formulation. Internal neuron variables are evaluated through the network mapping rather than introduced as independent optimization variables, while McCormick relaxations and convex and concave envelopes of the activation functions support deterministic global optimization. Their implementation in MAiNGO showed computational advantages over BARON on the reported examples. While reduced-space formulations address computational complexity by reducing the number of explicit optimization variables, piecewise-linear approximations address it by replacing nonlinear activation equations with linear segments. The latter approach underlies the PWL-tanh formulations in PWLNet and enables mixed-integer linear optimization when the remaining objective and constraints are linear.}

ReLU neural networks are particularly attractive for optimization because the ReLU activation is piecewise affine. Grimstad and Andersson showed that trained ReLU networks can be embedded as surrogate models in mixed-integer linear programs using big-M formulations \citep{grimstad2019relu_mip}. In these formulations, each uncertain ReLU activation is represented using linear inequalities and a binary variable. This enables deterministic optimization over trained ReLU networks using standard MILP solvers. Deterministic global optimization methods and open-source software such as \texttt{reluMIP} further support the deployment of embedded ReLU networks \citep{schweidtmann2019embedded,lueg2021relumip}. Subsequent work has developed stronger formulations, partition-based formulations, and bound-tightening approaches to improve the computational performance of ReLU-MIP models \citep{anderson2020strong_relu, tsay2021partition_relu}.

The practical performance of ReLU-MIP formulations depends strongly on network size and bound quality. Each additional ReLU neuron can introduce additional variables, constraints, and binary decisions. Furthermore, loose big-M constants weaken the continuous relaxation and increase branch-and-bound effort. For this reason, compact architectures, tight variable bounds, and bound-tightening procedures are essential for practical neural-network optimization. Similar formulation-strength considerations also appear in neural-network verification, where MILP and convex-relaxation methods are used to certify network properties over bounded input regions \citep{katz2017reluplex, tjeng2019verifying_relu}. Although the objective of verification differs from process optimization, the underlying mathematical challenges are closely related.

Open-source algebraic modeling environments have helped make neural-network embedding more accessible. Pyomo provides a flexible Python-based modeling framework for mathematical programming \citep{hart2011pyomo, bynum2021pyomo}. Building on Pyomo, OMLT provides tools for embedding trained machine-learning models, including neural networks and gradient-boosted trees, into larger Pyomo optimization problems \citep{ceccon2022omlt}. These tools substantially reduce the implementation burden associated with neural-network embedding. However, full process-surrogate workflows still require additional components such as data generation, scaling, formulation selection, pruning, bound management, solver selection, and validation of optimized decisions on the original process model.

\subsection{Formulation Strength, Bounds, and Solver Choice}

The same trained neural network can lead to different optimization problem classes depending on the activation function, approximation strategy, and objective function. ReLU networks with linear objectives and linear process constraints typically lead to MILP formulations. Smooth tanh networks can be embedded exactly as nonlinear programming formulations, but these NLPs are generally nonconvex and may have many local optima. Piecewise-linear approximations of tanh activations can instead be represented using MILP, SOS2, or SOSX-style formulations \citep{keha2004pwl, koksal2026sosx_training}. If the objective contains squared tracking errors, as commonly occurs in control applications, the resulting problem may become a mixed-integer quadratic program or may be reformulated as a mixed-integer second-order cone program.

This solver-dependent behavior is important in process systems engineering because computational performance is often as important as prediction accuracy. A neural network with excellent test-set accuracy may be impractical if its embedded optimization formulation is too large or numerically weak. Conversely, a slightly less accurate but much smaller surrogate may produce better overall decision performance if it enables faster and more reliable optimization. Therefore, formulation choice should be viewed as part of the surrogate-modeling workflow rather than as a post-processing detail.

Bounds and scaling are particularly important in this context. Neural networks are commonly trained on scaled data, but optimization models are usually written in physical units. The scaling transformations must therefore be represented consistently when the network is embedded into an optimization model. In addition, valid lower and upper bounds are required for big-M formulations, interval propagation, and piecewise-linear activation approximations. Tight bounds improve formulation strength, reduce relaxation gaps, and can remove inactive neurons or unused activation regions \citep{grimstad2019relu_mip}. These considerations motivate bound-aware and optimization-aware software workflows for neural surrogates.

\subsection{Piecewise-Linear tanh Networks and SOS2/SOSX Methods}

Although ReLU networks are convenient for MILP embedding, they are not always the most compact representation of chemical process behavior. Many process responses are smooth, saturating, and sigmoidal due to reaction kinetics, equilibrium effects, limiting reactants, and transport phenomena. Such behavior can often be represented more naturally by tanh activations than by ReLU activations. A sufficiently large ReLU network can approximate the same behavior, but doing so may require more neurons and potentially deeper architectures. This can increase the size of the corresponding optimization formulation.

Piecewise-linear tanh networks provide an alternative strategy. Instead of embedding the smooth tanh activation exactly, the activation is approximated over a bounded region using a set of linear segments. Each activation can then be represented using convex-combination variables, segment-selection variables, SOS2 constraints, or related mixed-integer formulations. 

Related MILP-based approaches have also combined piecewise-linear approximations with neural-network training and feature selection \citep{sildir2022milp_training}. The number and placement of breakpoints determine the trade-off between approximation accuracy and optimization model size. More breakpoints improve the approximation of the smooth activation, but they also increase the number of variables and constraints in the embedded model.

The development of piecewise-linear neural-network formulations for process optimization has been advanced by the work of K\"oksal, Ayd{\i}n, and T\"urkay. K\"oksal and Ayd{\i}n proposed physics-informed piecewise-linear neural networks for process optimization, demonstrating the value of combining process knowledge with optimization-compatible neural surrogates \citep{koksal2023physics_pwl}. K\"oksal, T\"urkay, and Ayd{\i}n later introduced a computationally efficient piecewise-linear training algorithm for neural networks using continuous special ordered sets \citep{koksal2026sosx_training}. These contributions provide the methodological foundation for training and optimizing PWL-tanh networks in a deterministic optimization setting.

Special ordered sets are central to these formulations. In SOS2 formulations, a piecewise-linear function is represented as a convex combination of adjacent breakpoints. The adjacency condition can be imposed through solver-native SOS2 constraints or equivalent mixed-integer formulations \citep{keha2004pwl, keha2006branch_cut_pwl}. SOSX-style approaches enforce the same validity condition through branching on continuous convex-combination variables rather than by introducing all binary variables explicitly. This is particularly useful in restricted training settings, such as bounded, fixed-readout, single-hidden-layer PWL-tanh networks trained under mean absolute error, where each branch-and-bound node can be solved as a linear program.

\subsection{Compact and Physics-Informed Neural Surrogates}

Network size is a central determinant of optimization tractability. A network that is accurate for prediction may still be unnecessarily large for optimization. Each additional neuron increases the number of variables and constraints in the embedded formulation, and in mixed-integer formulations it may also increase the number of binary variables. Therefore, neural-network compression is not only useful for reducing inference cost; it is also directly relevant to mathematical-programming tractability.

Classical pruning methods remove weights or neurons that contribute little to model accuracy. Optimal Brain Damage introduced a second-order saliency criterion for removing parameters while limiting the increase in training error \citep{lecun1989obd}. Later work developed related second-order pruning methods, practical Optimal Brain Damage variants, magnitude-based pruning, structured pruning, and retraining-based compression pipelines \citep{hassibi1993obs,liu2014obd,han2015deep_compression,cheng2023pruning_survey}. In optimization-aware surrogate modeling, pruning has an additional interpretation: it reduces the algebraic size of the optimization formulation. Dead ReLU removal, always-active ReLU simplification, saturated tanh pruning, magnitude pruning, neuron pruning, and group-regularized compression can all reduce the downstream optimization burden.

Recent work has examined more directly how pruning-induced sparsity changes the
computational behavior of optimization problems containing embedded ReLU
networks. Through comprehensive empirical studies across multiple network
architectures, an illustrative test function, and a chemical-process flowsheet,
\citet{lastrucci2026pruning} show that weight and node pruning can accelerate
deterministic optimization through three complementary mechanisms: reducing the
algebraic problem size, decreasing the number of integer variables, and
tightening the pre-activation bounds used in big-$M$ formulations. Their results
also show that the most effective pruning strategy depends on network
architecture: weight pruning is particularly effective for deep, narrow
networks, whereas node pruning is more effective for shallow, wider networks.
Based on these computational results, they recommend treating pruning as a
standard step when developing neural-network surrogates for embedded
optimization, particularly in engineering applications that require repeated
solves.

Physics-informed learning provides another mechanism for improving surrogate usefulness in process optimization. In process systems engineering, partial physical knowledge is often available even when a complete first-principles model is not. This information may include mass balances, conservation relationships, operating bounds, monotonicity constraints, or known qualitative behavior. Incorporating such information into the loss function can reduce physically inconsistent predictions and improve extrapolation behavior in data-sparse regions. This is especially important in optimization, because the optimizer may drive the surrogate toward regions of the input space that are poorly represented in the training data. Physics-informed process surrogates are therefore evaluated not only by prediction error, but also by physical consistency and downstream decision quality \citep{koksal2023physics_pwl}.

\subsection{Dynamic Learned Surrogates for NMPC}

Dynamic process optimization and model predictive control introduce additional challenges for neural surrogates. In nonlinear model predictive control, a model is used to predict the future evolution of the process over a finite horizon, and an optimization problem is solved repeatedly in a receding-horizon fashion. Learned surrogates can be used in this setting when first-principles dynamic models are unavailable, expensive, or difficult to embed directly.

Several neural architectures have been used for dynamic process modeling. Feedforward transition models approximate one-step or multi-step maps of the form
\begin{equation}
x_{k+1} = f_{\theta}(x_k,u_k),
\end{equation}
where $x_k$ is the current state and $u_k$ is the manipulated input. NARX models use lagged inputs and outputs to represent dynamic behavior without explicitly introducing a recurrent hidden state. Recurrent neural networks introduce an internal hidden state that evolves over time, making them suitable for sequence modeling. Long short-term memory networks were developed to improve learning over long time horizons by addressing the vanishing-gradient limitations of standard recurrent networks \citep{hochreiter1997lstm}. GRUs and other gated recurrent models provide related sequence-modeling capabilities.

Recent reviews and studies have shown increasing interest in neural-network and machine-learning models for MPC \citep{ren2022nn_mpc_review, wu2024ml_mpc_review}. However, embedding dynamic neural surrogates into NMPC problems can be computationally expensive because the network must be unrolled over the prediction horizon. The size of the optimization problem therefore grows with both the network architecture and the control horizon. In ReLU-based MPC, this can lead to a large MILP; in smooth recurrent models, it can lead to large nonconvex NLPs. These issues motivate compact neural architectures, pruning, bound tightening, piecewise-linear formulations, and careful selection of the optimization backend.

Learned NMPC also raises control-specific concerns. Prediction accuracy over a test set does not guarantee closed-loop performance. Constraint satisfaction, recursive feasibility, stability, robustness to model mismatch, and extrapolation outside the training data are all important considerations. PWLNet does not attempt to provide general closed-loop stability guarantees for arbitrary learned surrogates. Instead, it focuses on producing compact optimization-compatible surrogate formulations and validating the resulting decisions against process models or simulators.

\subsection{Convex Neural-Network Architectures}

Another related direction is the development of neural-network architectures that are convex by construction. Input-convex neural networks impose structural restrictions on the network, including nonnegative weights in selected layers and convex nondecreasing activation functions, so that the output is convex with respect to specified inputs \citep{amos2017icnn}. Such architectures are attractive when the learned function must be optimized efficiently, because convexity can be preserved in the downstream optimization problem.

However, convex neural-network architectures restrict the class of functions that can be represented. Many chemical process optimization problems are nonconvex due to reaction networks, recycle structures, phase behavior, competing pathways, and dynamic effects. For such systems, enforcing global convexity may be too restrictive. General ReLU and tanh networks provide broader representational flexibility, but require mixed-integer, nonlinear, conic, or piecewise-linear reformulations for optimization. PWLNet is positioned in this more general setting: it supports flexible neural surrogates while providing tools to make the resulting optimization problems more tractable.

\subsection{Software Gap and Positioning of PWLNet}

Several existing tools support parts of the optimization-aware neural-surrogate workflow. TensorFlow and Keras provide neural-network training \citep{abadi2016tensorflow,chollet2015keras}. Pyomo provides algebraic optimization modeling \citep{hart2011pyomo}. OMLT enables trained machine-learning models to be embedded in Pyomo optimization problems \citep{ceccon2022omlt}. ReLU-MIP formulations provide a route for deterministic optimization with ReLU neural networks \citep{grimstad2019relu_mip}. PWL and SOS-based methods provide a route for training and optimizing piecewise-linear neural networks \citep{keha2006branch_cut_pwl,koksal2026sosx_training}. However, these components are often used through case-specific scripts, making it difficult to compare activation functions, formulations, pruning strategies, solvers, and validation metrics under a unified workflow.

This creates a reproducibility and usability gap. Many published surrogate-optimization studies demonstrate a specific model, formulation, or case study, but do not provide a reusable framework for moving from training data to embedded optimization and exact-process validation. As a result, researchers often need to rebuild scaling, model export, Pyomo embedding, solver configuration, pruning, and benchmarking pipelines for each new process system.

PWLNet addresses this gap by providing an open-source Python framework for optimization-compatible neural surrogate modeling in process systems engineering. The framework builds on ReLU-MIP surrogate optimization, Pyomo/OMLT neural-network embedding, optimization-aware pruning, physics-informed process surrogates, and the PWL/SOSX methodology of K\"oksal, Ayd{\i}n, and T\"urkay \citep{koksal2026sosx_training}. Rather than promoting a single activation function or formulation, PWLNet emphasizes interchangeability. Users can train ReLU or tanh networks, embed exact or piecewise-linear formulations, adjust breakpoint placement and formulation size, apply pruning and compression, select among MILP, MIQP, MISOCP, NLP, and SOSX-style workflows, and validate optimized decisions on process case studies. In this way, PWLNet connects neural networks as predictive models with neural networks as reusable, solver-compatible components of optimization and control workflows.

\section{PWLNet Framework}
\label{sec:framework}

PWLNet is organized as a layered Python package for constructing neural
surrogates that remain usable after training as algebraic optimization
models. PWLNet is designed around a clear separation between training, representation, and optimization: neural networks may be trained using different machine-learning frameworks, stored in a framework-independent internal format, and then converted into different optimization formulations for different solver technologies.  A
trained model can therefore move through the workflow
\emph{train--scale--bound--compress--formulate--solve} without retaining a
dependency on the framework used to estimate its parameters.  Dependencies
on Keras, TensorFlow, Pyomo, OMLT, and commercial solvers are imported only by
the modules that need them.  The base package consequently remains lightweight
and the licensed solver interfaces remain optional.

PWLNet is organized as a modular workflow in which each stage can be varied independently to test different neural architectures, compression choices, reformulation strategies, and solver configurations on the same benchmark problem:
\begin{equation}
\text{benchmark data}
\rightarrow \text{surrogate training}
\rightarrow \text{compression}
\rightarrow \text{reformulation}
\rightarrow \text{pyomo embedding}
\rightarrow \text{optimization}.
\end{equation}

\subsection{Facade and Framework-Independent Representation}

The user-facing \texttt{PWLNet} class is a thin facade rather than a second
model implementation. It delegates work through a wrapped \texttt{NetworkDefinition}, including training, prediction,serialization, summaries, bound propagation, pruning, compression, and
algebraic export. Class methods expose dense Keras training, physics-informed training, SOSX training,
and recurrent training, while instance methods such as \texttt{predict},
\texttt{save}, \texttt{bounds}, \texttt{compress}, and \texttt{to\_pyomo}
provide the common downstream workflow.  The fitted trainer is retained when
appropriate so that convergence histories and training diagnostics remain
available, but the wrapped network is the authoritative deployable object.

The framework-independent \texttt{NetworkDefinition} is implemented using
NumPy \citep{harris2020numpy} and contains an ordered sequence of \texttt{DenseLayer} objects.  Each
layer stores a weight matrix, a bias vector, and an activation identifier
(ReLU, hyperbolic tangent, or linear).  For a scaled input
$\hat{\boldsymbol{x}}$, layer $k$ evaluates
\begin{equation}
  \boldsymbol{h}_{k+1}
  = \phi_k\!\left(\boldsymbol{h}_k\boldsymbol{W}_k
  + \boldsymbol{b}_k\right),
  \qquad \boldsymbol{h}_0=\hat{\boldsymbol{x}}.
\end{equation}
The object validates layer connectivity, scaler dimensions, and input bounds;
supports prediction in physical units; and provides a JSON-compatible,
versioned serialization of layers, scalers, and bounds.  It deliberately has
no imports from Keras, TensorFlow, Pyomo, or OMLT.  This boundary makes the
representation portable, testable, and suitable as the shared contract among
training, compression, and optimization backends.

\subsection{Scaling and Bound Propagation}

Input and output transformations are represented explicitly by
\texttt{OffsetScaler}.  For each feature, PWLNet uses the affine convention
\begin{equation}
  \hat{x}_j = \frac{x_j-o_j}{s_j},
  \qquad
  x_j = s_j\hat{x}_j+o_j,
\end{equation}
which covers standardization, min--max scaling, and the identity mapping.  The
same convention is used by OMLT \citep{ceccon2022omlt}, allowing PWLNet to expose optimization-block
inputs and outputs in original engineering units while formulating the neural
network internally in scaled coordinates.  The scalers are fitted only on the
training partition and are stored with the exported network.

Finite input bounds are also part of \texttt{NetworkDefinition} and are stored
in original units.  The standard dense trainer initializes them from the
componentwise extrema of the training inputs; an application may replace them
with a validated operating envelope before optimization.  The
\texttt{propagate\_bounds} routine first transforms this box into scaled space
and then applies interval arithmetic to every affine layer.  Splitting the
weight matrix into positive and negative parts gives the pre-activation bounds
\begin{align}
  \underline{\boldsymbol{z}}
  &=\underline{\boldsymbol{h}}\boldsymbol{W}^{+}
    +\overline{\boldsymbol{h}}\boldsymbol{W}^{-}+\boldsymbol{b},\\
  \overline{\boldsymbol{z}}
  &=\overline{\boldsymbol{h}}\boldsymbol{W}^{+}
    +\underline{\boldsymbol{h}}\boldsymbol{W}^{-}+\boldsymbol{b}.
\end{align}
Because the supported activations are monotonic, their interval images follow
from evaluating the two endpoints.  The resulting bounds supply output-range
estimates, big-$M$ constants, breakpoint domains, and certificates for exact
ReLU simplifications.  They are exact for an individual affine map over a box
but can be conservative after several layers.

For safety, the input bounds used in an optimization study should represent the
validated process operating region rather than an invitation to extrapolate.
The serialized \texttt{NetworkDefinition}, explicit scaling convention, and
backend-independent solver interface make that domain assumption inspectable
and reproducible.

\subsection{Training Modules}

Training is defined by a family of adapters that all export optimization-ready
objects.  \texttt{KerasTrainer} fits dense ReLU or tanh networks in scaled
space, extracts the learned parameters into NumPy \texttt{DenseLayer}
instances, and checks numerical equivalence between the Keras model and the
exported \texttt{NetworkDefinition}.  \texttt{PhysicsInformedTrainer} extends
the supervised objective with differentiable, process-specific residuals.
Physics terms can be supplied through a \texttt{PhysicsLoss} implementation,
a declarative \texttt{PhysicsSpec}, or a registered provider, thereby keeping
the generic optimization loop independent of process equations.

\hl{The training workflow also supports PWL-native networks, in which the
piecewise-linear tanh activation is used during fitting and retained during
prediction and optimization embedding. This avoids the activation-replacement
error introduced when a smoothly trained tanh network is linearized afterward.
In the hydroformylation benchmark, PWL-native fitting is implemented in the
study-specific training workflow rather than as a general option shared by the
trainer adapters. The distinction between these deployment paths is detailed
in Section}~\ref{subsec:pwl_tanh}.

\texttt{SosxTrainer} provides an alternative global-training path for a
bounded, single-hidden-layer, piecewise-linear tanh network.   It minimizes mean absolute error with a fixed linear readout, admitting an exact branch-and-bound solution over linear-program subproblems. Its diagnostics include the incumbent, optimality
gap, node count, and solver and wall-clock times.  Finally,
\texttt{RNNTrainer} estimates a ReLU Elman cell and exports a framework-free
\texttt{RecurrentCell}; the cell can be rolled out numerically or unrolled over
a control horizon as an exact mixed-integer linear model.  Thus, feed-forward,
physics-informed, globally trained PWL, and recurrent models share the same
architectural principle: framework-specific estimation followed by
framework-independent deployment.

The benchmark training scripts set the NumPy and TensorFlow random seeds and
enable TensorFlow operation determinism, so rerunning a recorded seed produces
identical trained weights on the supported software stack. The serialized
networks used by the executable notebooks are stored with their dataset and
training settings. The transition datasets are likewise committed together
with their generation settings, providing a direct record from simulated data
through training and optimization.

\subsection{Pruning and Compression}

Compression operates on the exported representation and never mutates the
source network.  The low-level pruning package contains structural operations
for copying a network, removing neurons, fusing adjacent linear layers, and
zeroing weights.  Bound-based ReLU routines remove neurons that are inactive
throughout the input box and replace always-active ReLUs by their affine maps;
these transformations preserve the represented function over the declared
domain.  Additional methods fold saturated tanh units, prune small-magnitude
weights or low-importance neurons, use group-lasso structure, or apply Optimal
Brain Damage saliencies.  Keras-based fine-tuning can subsequently recover
predictive accuracy.

The facade's \texttt{compress} method composes these operations into a staged
pipeline: exact ReLU simplification, optional weight pruning, optional neuron
pruning, and optional fine-tuning.  Each stage returns a
\texttt{PruneOutcome}, while reporting utilities compare parameter counts,
prediction error, and estimated formulation size before and after compression.
This arrangement makes the effect of compression visible both in statistical
accuracy and in the size of the downstream optimization problem.

\subsection{Formulation Backends and Solver Interfaces}

Formulation code is partitioned mathematically. The OMLT
backend converts \texttt{NetworkDefinition} directly, layer by layer, without
a Keras file-format round trip.  It provides a ReLU big-$M$ formulation for
mixed-integer linear optimization and an exact smooth tanh formulation for
nonlinear programming.  Scaling objects and scaled input bounds are transferred
to OMLT so that the resulting \texttt{block.inputs} and
\texttt{block.outputs} retain physical units.

The custom Pyomo-PWL backend represents tanh activations using either SOS2 or
big-$M$ piecewise-linear constraints, with reusable breakpoint and epigraph
utilities.  The SOSX backend separates the SOS2-dropped node relaxation,
violation detection, breakpoint-window branching, node selection, incumbent
management, and fathoming logic.  This supports branch-and-LP or
branch-and-MILP solution through \texttt{solve\_sosx}, including warm starts
for repeated problems such as model predictive control.  Recurrent models use
the same formulation layer after their cells are unrolled over the prediction
horizon.

Solver selection is isolated behind
\texttt{get\_\allowbreak pyomo\_\allowbreak solver},
\texttt{get\_\allowbreak milp\_\allowbreak solver},
\texttt{get\_\allowbreak conic\_\allowbreak solver}, and
\texttt{get\_\allowbreak nlp\_\allowbreak solver}.  These functions register Pyomo plugins, check
availability, apply options across both classic and APPSI interfaces, and
produce actionable errors when an optional executable or license is missing.
HiGHS through \texttt{appsi\_highs} is the default open-source LP/MILP
interface \citep{huangfu2018}.  Gurobi and MOSEK remain optional licensed alternatives, and Ipopt
is discovered only when an exact nonlinear formulation is requested \citep{gurobi2026,mosek2026,wachter2006}.  This
layer prevents solver-specific setup from leaking into the network and process
modules.

\subsection{Process Modules and Executable Notebooks}

Process-specific equations are organized below \texttt{pwlnet.processes} and
are kept distinct from the reusable learning and formulation layers.  The
Williams--Otto package contains the steady-state reactor equations, parameter
data, simulation and dataset generation, profit calculation, physics-informed
loss terms, exact-reference optimization, and surrogate-optimization models.
The hydroformylation package contains a dynamic semibatch model, numerical
integration, transition-data generation, delayed-measurement simulation,
open-loop optimal control, and shrinking- or receding-horizon NMPC helpers.
The Van de Vusse package similarly provides steady-state and dynamic reactor
models, data generators, direct and multiple-shooting control models, NMPC
utilities, and the Sahlodin benchmark variant.  These packages serve both as
reproducible benchmarks and as examples of the boundary between generic
surrogate infrastructure and application equations.

Fourteen notebooks provide executable, progressively richer examples.  The
first three compare ReLU/OMLT, exact tanh/OMLT, and PWL-tanh formulations on the
Williams--Otto problem.  Notebook 4 develops hydroformylation PWL-NMPC, and
Notebook 5 studies ReLU pruning and its effects on accuracy, model size, and
solution time.  Notebooks 6 and 7 demonstrate physics-informed training for
Williams--Otto and hydroformylation, respectively.  Notebooks 8--10 cover the
SOSX branch-and-LP workflow, short and full-horizon hydroformylation control, and a
Williams--Otto input-noise stress test.  Notebooks 11 and 12 extend the control
examples to Van de Vusse NMPC using, respectively, feed-forward PWL surrogates
and a pruned, physics-informed ReLU recurrent network embedded in an MPC MILP. Notebook 13 applies the complete workflow to the Boston housing regression benchmark, and notebook 14 adds the Van de Vusse policy-validation studies.
Together, the notebooks function as integration tests and as reproducible
records of the train--formulate--solve--validate workflow.

\subsection{Installation and Minimal Workflow}

The intended installation command for the released package is

\begin{verbatim}
pip install pwlnet
\end{verbatim}

The base installation contains NumPy/SciPy data structures and utilities
\citep{harris2020numpy,virtanen2020scipy}. The experimental workflow uses
pandas for tabular data, scikit-learn for preprocessing and evaluation, and
Matplotlib for visualization \citep{mckinney2010pandas,pedregosa2011scikit,hunter2007matplotlib}.
Training and optimization dependencies are optional; the complete open-source
workflow used below can be installed with

\begin{verbatim}
pip install "pwlnet[training,highs]"
\end{verbatim}

The following minimal example generates a two-input regression dataset,
trains and serializes a ReLU surrogate, exports it as an OMLT/Pyomo big-$M$
block, and maximizes its predicted output with HiGHS.  Both the optimization
variables and the reported solution use the original, unscaled units.

\begin{verbatim}
import numpy as np
import pyomo.environ as pyo

from pwlnet import PWLNet
from pwlnet.solvers import get_pyomo_solver

# Training data in original units.
rng = np.random.default_rng(7)
X = rng.uniform([-1.0, -2.0], [1.0, 2.0], size=(200, 2))
y = (1.5 * np.maximum(X[:, 0], 0.0) - 0.25 * X[:, 1])[:, None]

# Train and export a framework-independent JSON representation.
net = PWLNet.train(
    X,
    y,
    hidden_layers=[8],
    activation="relu",
    epochs=200,
    batch_size=32,
    validation_split=0.1,
    random_seed=7,
    verbose=0,
)
net.save("surrogate.json")

# Export the trained network to an optimization model.
model, surrogate = net.to_pyomo(
    backend="omlt",
    formulation="relu_bigm",
    block_name="surrogate",
)
model.objective = pyo.Objective(
    expr=surrogate.outputs[0],
    sense=pyo.maximize,
)

# Solve over the input box stored in NetworkDefinition.
solver = get_pyomo_solver("appsi_highs")
results = solver.solve(model)
pyo.assert_optimal_termination(results)

x_star = np.array([
    pyo.value(surrogate.inputs[i]) for i in surrogate.inputs
])
y_star = pyo.value(surrogate.outputs[0])
print("x* =", x_star)
print("surrogate y* =", y_star)
\end{verbatim}

\section{Methods}
\subsection{Linearized Neural Surrogate Formulations}
\label{sec:formulations}

This section presents the mathematical formulations supported by
PWLNet: the ReLU big-$M$ mixed-integer formulation for feedforward and
recurrent networks, the exact smooth tanh nonlinear formulation, the
piecewise-linear tanh big-$M$ and SOS2 formulations, the SOSX
branch-and-bound approach over continuous convex-combination
variables, and the quadratic and second-order-cone objective variants
used in control applications.

Consider a feedforward neural network surrogate
\begin{equation}
    \hat{y} = f_L \left( W_L f_{L-1} \left( \cdots f_1(W_1 x + b_1) \right) + b_L \right),
    \label{eq:nn_general}
\end{equation}
where $x$ is the input vector, $\hat{y}$ is the network prediction, $W_\ell$ and $b_\ell$ are the weights and biases of layer $\ell$, and $f_\ell$ is the activation function.

PWLNet focuses on reformulations that make Eq.~\eqref{eq:nn_general} compatible with algebraic optimization models. The main formulations considered in this study are summarized in Table~\ref{tab:formulations}. If squared-error objectives are used within optimization-based training or model fitting, they are treated as quadratic objectives, leading to QP or MIQP formulations depending on the presence of integer variables; the same convex quadratic objective can equivalently be reformulated with a second-order cone constraint, yielding MISOCP models for conic solvers such as MOSEK (Section~\ref{sec:PWL-mse}).

\subsection{ReLU Mixed-Integer Linear Formulation}
\label{subsec:relu_mip}

For networks with rectified linear unit (ReLU) activations, each neuron can be reformulated exactly using mixed-integer linear constraints, provided that valid bounds are available for its pre-activation value. For a scalar ReLU activation,
\begin{equation}
z = \max(0,s),
\label{eq:relu}
\end{equation}
where $s$ denotes the pre-activation value and $z$ denotes the post-activation output. Given valid bounds $L \leq s \leq U$, the ReLU relation can be represented by the standard big-$M$ formulation
\begin{equation}
\begin{aligned}
z &\geq s, \\
z &\geq 0, \\
z &\leq s - L(1-\delta), \\
z &\leq U\delta, \\
\delta &\in \{0,1\}.
\end{aligned}
\label{eq:relu_bigm}
\end{equation}
Here, the binary variable $\delta$ indicates whether the neuron is active. When $\delta = 1$, the formulation enforces $z=s$ with $s \geq 0$; when $\delta = 0$, it enforces $z=0$ with $s \leq 0$. Applying this formulation to every ReLU neuron converts a trained ReLU neural network into an equivalent mixed-integer linear representation. In PWLNet, this reformulation is implemented primarily through the OMLT/Pyomo backend.

\begin{table}[h]
\centering
\caption{Neural surrogate formulations considered in PWLNet.}
\label{tab:formulations}
\begin{tabular}{p{0.18\textwidth} p{0.30\textwidth} p{0.22\textwidth} p{0.18\textwidth}}
\toprule
\textbf{Activation} & \textbf{Formulation} & \textbf{Backend} & \textbf{Problem class} \\
\midrule
ReLU & big-$M$ or partition formulation & OMLT/Pyomo & MILP \\
ReLU & quadratic prediction loss & Pyomo & MIQP \\
tanh & exact nonlinear activation & OMLT/Pyomo & NLP \\
tanh & piecewise-linear big-$M$ & Pyomo & MILP or MIQP \\
tanh & piecewise-linear SOS2 & Pyomo & MILP or MIQP \\
tanh & piecewise-linear big-$M$/SOS2 with MSE objective & Pyomo/MOSEK & MIQP or MISOCP \\
piecewise tanh & SOSX & custom brancher/Pyomo & branch-and-LP / branch-and-MILP \\
\bottomrule
\end{tabular}
\end{table}

\subsection{ReLU Recurrent Neural Network Mixed-Integer Linear Formulation}
\label{subsec:relu_rnn_mip}

PWLNet also supports recurrent neural network surrogates for dynamic prediction and control. In particular, a ReLU Elman recurrent neural network (RNN) can be embedded into an optimization model by unrolling the recurrent cell over a finite prediction horizon and applying the ReLU big-$M$ formulation at each time step.

Let $u_t \in \mathbb{R}^{n_u}$ denote the exogenous input or control vector at time step $t$, $h_t \in \mathbb{R}^{n_h}$ denote the hidden state, and $y_t \in \mathbb{R}^{n_y}$ denote the observable output. A single ReLU Elman cell is given by
\begin{equation}
h_{t+1}
=
\operatorname{ReLU}
\left(
u_t^\top W_{xh}
+
h_t^\top W_{hh}
+
b_h
\right),
\label{eq:relu_rnn_cell}
\end{equation}
where $W_{xh} \in \mathbb{R}^{n_u \times n_h}$ is the input-to-hidden weight matrix, $W_{hh} \in \mathbb{R}^{n_h \times n_h}$ is the hidden-to-hidden recurrent weight matrix, and $b_h \in \mathbb{R}^{n_h}$ is the hidden bias vector. The output is computed through a linear readout,
\begin{equation}
y_t = h_{t+1}^\top W_{hy} + b_y,
\label{eq:relu_rnn_readout}
\end{equation}
where $W_{hy} \in \mathbb{R}^{n_h \times n_y}$ and $b_y \in \mathbb{R}^{n_y}$.

The recurrent update in Eq.~\eqref{eq:relu_rnn_cell} can be written as a standard dense ReLU layer by concatenating the current input and hidden state. Defining
\begin{equation}
\xi_t =
\begin{bmatrix}
u_t \
h_t
\end{bmatrix},
\qquad
\widetilde{W}
=
\begin{bmatrix}
W_{xh} \\
W_{hh}
\end{bmatrix},
\label{eq:rnn_stacking}
\end{equation}
the hidden-state update becomes
\begin{equation}
h_{t+1}
=
\operatorname{ReLU}
\left(
\xi_t^\top \widetilde{W}
+
b_h
\right).
\label{eq:rnn_stacked_relu}
\end{equation}
This representation is useful because each recurrent step has the same mathematical form as a feedforward dense ReLU layer. Therefore, the same big-$M$ formulation used for feedforward ReLU networks can be applied directly to each hidden unit at each time step.

For hidden unit $j$ at time step $t$, let $s_{t,j}$ denote the pre-activation value and $h_{t+1,j}$ denote the corresponding post-activation value. The affine pre-activation is
\begin{equation}
s_{t,j}
=
\sum_{i=1}^{n_u}
W_{xh,ij} u_{t,i}
+
\sum_{k=1}^{n_h}
W_{hh,kj} h_{t,k}
+
b_{h,j}.
\label{eq:rnn_preactivation}
\end{equation}
Given valid bounds $L_{t,j} \leq s_{t,j} \leq U_{t,j}$, the ReLU relation
$h_{t+1,j} = \max(0,s_{t,j})$ is represented by
\begin{equation}
\begin{aligned}
h_{t+1,j} &\geq s_{t,j}, \\
h_{t+1,j} &\geq 0, \\
h_{t+1,j} &\leq s_{t,j} - L_{t,j}(1-\delta_{t,j}), \\
h_{t+1,j} &\leq U_{t,j}\delta_{t,j}, \\
\delta_{t,j} &\in \{0,1\}.
\end{aligned}
\label{eq:rnn_bigm}
\end{equation}
The binary variable $\delta_{t,j}$ indicates whether hidden unit $j$ is active at time step $t$. If $\delta_{t,j}=1$, the formulation enforces $h_{t+1,j}=s_{t,j}$ with $s_{t,j}\geq 0$. If $\delta_{t,j}=0$, it enforces $h_{t+1,j}=0$ with $s_{t,j}\leq 0$.

To embed the RNN over a prediction horizon $H$, the recurrent cell is unrolled for $t=0,\ldots,H-1$. The initial hidden state is fixed from the available past trajectory,
\begin{equation}
h_0 = h_0^\ast,
\label{eq:rnn_initial_hidden}
\end{equation}
and the optimization model includes the recurrent constraints
\begin{equation}
\begin{aligned}
s_{t,j}
&=
\sum_{i=1}^{n_u}
W_{xh,ij}u_{t,i}
+
\sum_{k=1}^{n_h}
W_{hh,kj}h_{t,k}
+
b_{h,j},
&& t=0,\ldots,H-1,; j=1,\ldots,n_h, \\
h_{t+1,j}
&=
\max(0,s_{t,j}),
&& t=0,\ldots,H-1,; j=1,\ldots,n_h, \\
y_{t,k}
&=
\sum_{j=1}^{n_h}
W_{hy,jk}h_{t+1,j}
+
b_{y,k},
&& t=0,\ldots,H-1,; k=1,\ldots,n_y,
\end{aligned}
\label{eq:rnn_unrolled}
\end{equation}
where each ReLU equality is implemented using the big-$M$ constraints above.

This unrolled formulation produces a mixed-integer linear representation of the trained ReLU RNN. The recurrent connection itself is linear in the previous hidden state; therefore, the only source of integer variables is the ReLU activation. In the worst case, the formulation introduces $H n_h$ binary variables, one for each hidden unit at each time step. This count can be reduced when bound propagation identifies hidden units that are always active or always inactive over the relevant operating region. If $L_{t,j}\geq 0$, the ReLU is always active and can be replaced by the linear equality $h_{t+1,j}=s_{t,j}$. If $U_{t,j}\leq 0$, the unit is always inactive and can be fixed to $h_{t+1,j}=0$.

The quality of the big-$M$ formulation depends on the tightness of the pre-activation bounds $L_{t,j}$ and $U_{t,j}$. In PWLNet, these bounds are obtained by recurrent interval bound propagation over the unrolled horizon. The resulting time-dependent bounds improve numerical conditioning, reduce unnecessary binary variables, and enable exact pruning of inactive ReLU units. Consequently, ReLU RNNs provide an optimization-compatible surrogate structure for dynamic prediction and model predictive control while preserving a mixed-integer linear problem class.

\subsection{Exact tanh Nonlinear Formulation}
\label{subsec:exact_tanh_nlp}

For neural networks with hyperbolic tangent activation functions, PWLNet also supports an exact nonlinear formulation. Unlike the ReLU and piecewise-linear tanh formulations, this approach does not introduce binary variables or approximate the activation function. Instead, the activation is represented directly through nonlinear equality constraints. For a tanh neuron,
\begin{equation}
z = \tanh(s),
\label{eq:tanh_activation}
\end{equation}
where $s$ is the pre-activation value and $z$ is the post-activation value.

\paragraph{Training formulation.}
In the training problem, the neural network weights and biases are decision variables. Given a dataset
$\{(x_i, y_i)\}_{i=1}^{N_s}$, the exact tanh training problem can be written as
\begin{equation}
\begin{aligned}
\min_{{W^\ell,b^\ell,s_i^\ell,z_i^\ell}}
\quad &
\frac{1}{N_s}\sum_{i=1}^{N_s}
\left| \hat{y}_i - y_i \right|_2^2
+ \mathcal{R}(W,b) \\
\text{s.t.}\quad
& z_i^0 = x_i,
&& i=1,\ldots,N_s, \\
& s_i^\ell = W^\ell z_i^{\ell-1} + b^\ell,
&& i=1,\ldots,N_s,;\ell=1,\ldots,L, \\
& z_i^\ell = \tanh(s_i^\ell),
&& i=1,\ldots,N_s,;\ell \in \mathcal{H}, \\
& z_i^\ell = s_i^\ell,
&& i=1,\ldots,N_s,;\ell \notin \mathcal{H}, \\
& \hat{y}_i = z_i^L,
&& i=1,\ldots,N_s.
\end{aligned}
\label{eq:tanh_training_obj}
\end{equation}
Here, $\mathcal{H}$ denotes the set of hidden layers with tanh activations, while the final layer is typically linear. The term $\mathcal{R}(W,b)$ represents optional regularization. Because the tanh activation is smooth and nonlinear, the resulting training problem is a nonlinear programming (NLP) problem. In the exact nonlinear training path, this NLP can be solved using a local nonlinear optimizer such as IPOPT.

\paragraph{Optimization embedding formulation.}
After training, the weights and biases are fixed. The trained neural network is then embedded into a downstream process optimization problem as a nonlinear surrogate model. Let $x$ denote the process decision variables or surrogate inputs, and let $\hat{y}$ denote the neural-network prediction. The embedded network is represented as
\begin{equation}
\begin{aligned}
z^0 &= x, \\
s^\ell &= W^\ell z^{\ell-1} + b^\ell,
&& \ell=1,\ldots,L, \\
z^\ell &= \tanh(s^\ell),
&& \ell \in \mathcal{H}, \\
z^\ell &= s^\ell,
&& \ell \notin \mathcal{H}, \\
\hat{y} &= z^L.
\end{aligned}
\label{eq:tanh_embed}
\end{equation}
The downstream optimization problem can therefore be written as
\begin{equation}
\begin{aligned}
\min_{x,\hat{y},s,z}
\quad & F(x,\hat{y}) \\
\text{s.t.}\quad
& g(x,\hat{y}) \leq 0, \\
& z^0 = x, \\
& s^\ell = W^\ell z^{\ell-1} + b^\ell,
&& \ell=1,\ldots,L, \\
& z^\ell = \tanh(s^\ell),
&& \ell \in \mathcal{H}, \\
& z^\ell = s^\ell,
&& \ell \notin \mathcal{H}, \\
& \hat{y} = z^L.
\end{aligned}
\label{eq:tanh_opt}
\end{equation}
Here, $F(x,\hat{y})$ is the downstream process objective and $g(x,\hat{y})$ represents process constraints, operating limits, or product-quality specifications. Since the trained weights and biases are fixed, the neural network contributes only nonlinear equality constraints to the optimization model. The resulting embedded problem is an NLP when the remaining process model is continuous.

In PWLNet, the exact tanh embedding is implemented through the OMLT/Pyomo backend. The trained tanh network is exported as an algebraic Pyomo block, where each tanh activation is represented exactly rather than through a piecewise-linear approximation. The resulting nonlinear optimization model is solved using IPOPT.

\subsubsection{\hl{Motivation for Piecewise-Linear tanh Surrogates}}

Exact tanh embeddings preserve smooth activations but generally produce
nonconvex nonlinear optimization problems. Local NLP solvers such as IPOPT
do not generally certify global optimality, while deterministic global
solution can become computationally demanding as network size increases.
Reduced-space methods address this difficulty while retaining smooth
activations, as discussed in Section~\ref{sec:embedding}.

Piecewise-linear tanh provides a representation compatible with
mixed-integer linear optimization when the remaining objective and
constraints are linear. It retains the saturating shape of tanh through a
finite set of linear segments, with breakpoint number and placement
controlling representation flexibility and formulation size. Whether this
produces a more efficient model than a ReLU network or an exact tanh
embedding depends on the application, network architecture, and formulation.

PWLNet supports both post-training approximation of smooth tanh networks
and PWL-native training. These paths differ in whether the activation
changes after training, as detailed in Section~\ref{subsec:pwl_tanh}. A
global optimization certificate applies to the embedded surrogate problem;
validation against the original process remains necessary.

\subsection{Piecewise-Linear tanh Formulation}
\label{subsec:pwl_tanh}

\begin{figure}
\centering
\includegraphics[width=0.85\textwidth]{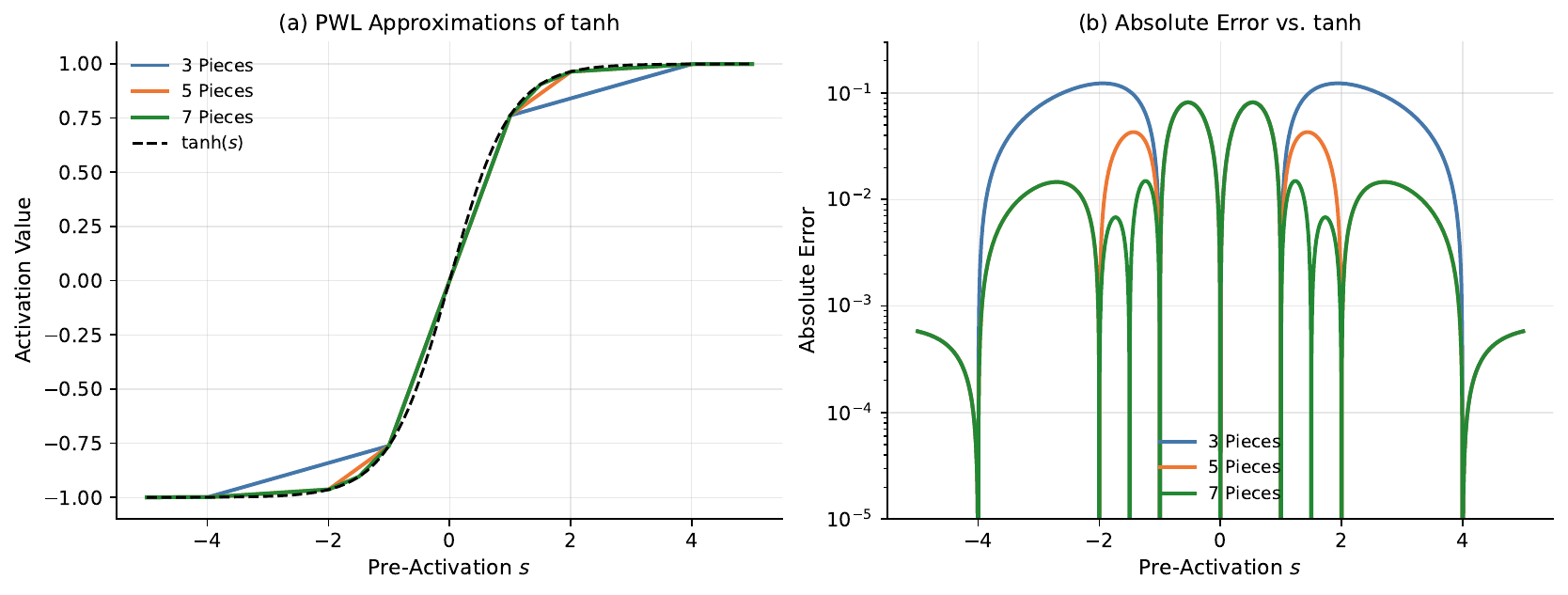}
\caption{Piecewise-linear approximations of the hyperbolic tangent activation function. The figure illustrates how a smooth tanh activation can be approximated by different numbers of linear segments before embedding the surrogate into an optimization model.}
\label{fig:tanh_pwl}
\end{figure}

\hl{The hyperbolic tangent activation function is useful for representing smooth nonlinear and saturating process behavior. However, its exact embedding introduces nonlinear activation equalities. PWLNet supports a piecewise-linear activation representation that can either replace smooth tanh after training or be used directly during training. The same piecewise-linear formulation applies to both cases, but only the latter preserves the activation representation consistently across training, prediction, and optimization embedding. We refer to this approach as PWL-native training.}

For a tanh neuron,
\begin{equation}
z = \tanh(s),
\end{equation}
where $s$ is the pre-activation value and $z$ is the post-activation value. The approximation is constructed on the interval
\begin{equation}
-4 \leq s \leq 4,
\end{equation}
because $\tanh(s)$ is already close to its saturation values outside this region. Thus, extending the approximation beyond this interval gives little additional accuracy while increasing the size of the optimization formulation.

Let the breakpoint vector be
\begin{equation}
a_0 < a_1 < \cdots < a_P,
\end{equation}
where $P$ is the number of linear pieces. For each segment $k=1,\ldots,P$, define
\begin{equation}
u_k = a_k - a_{k-1},
\qquad
g_k = \tanh(a_k)-\tanh(a_{k-1}).
\end{equation}
Then, for $s \in [a_{k-1},a_k]$, the piecewise-linear approximation is
\begin{equation}
\widehat{\tanh}(s)
=
\tanh(a_{k-1})
+
\frac{g_k}{u_k}
\left(s-a_{k-1}\right),
\qquad k=1,\ldots,P.
\label{eq:pwl_tanh_segment}
\end{equation}
Equivalently, each segment is represented by a linear function whose slope is determined by the corresponding neighboring breakpoint values.

In this work, the breakpoint sets are selected to preserve the structure of the tanh activation. Since tanh is symmetric around the origin, symmetric breakpoint sets with an odd number of linear pieces are used. This gives the same number of segments on the negative and positive sides of the origin, together with a central segment around zero. This choice preserves the qualitative shape of tanh and avoids asymmetric discretization artifacts that may occur when an even number of pieces is used.

The breakpoint sets considered in this study are
\begin{equation}
\begin{aligned}
3\text{-piece:}\quad
&[-4,-1,1,4], \\
5\text{-piece:}\quad
&[-4,-2,-1,1,2,4], \\
7\text{-piece:}\quad
&[-4,-2,-1.5,-1,1,1.5,2,4].
\end{aligned}
\label{eq:pwl_breakpoint_sets}
\end{equation}
For example, the 3-piece approximation is written explicitly as
\begin{equation}
\widehat{\tanh}_{3}(s) =
\begin{cases}
0.08s - 0.68, & -4 \leq s \leq -1, \\
0.76s,        & -1 \leq s \leq 1, \\
0.08s + 0.68, & 1 \leq s \leq 4.
\end{cases}
\label{eq:tanh_3_piece}
\end{equation}
The coefficients displayed in Eq.~\eqref{eq:tanh_3_piece} are rounded for readability; training, prediction, and optimization embedding use the same unrounded breakpoint values and segment coefficients.

The 5-piece and 7-piece approximations follow the same construction using the corresponding breakpoint vectors. Increasing the number of pieces improves the local approximation of tanh but also increases the number of linear regions that must be represented in the optimization model. Therefore, the number of pieces controls the trade-off between approximation accuracy and computational tractability which can be seen in Figure \ref{fig:tanh_pwl}. Further details on the piecewise-linear tanh approximation and its use in neural-network formulations can be seen in K\"oksal, T\"urkay, and Ayd{\i}n \citep{koksal2026sosx_training}.

\subsubsection{\hl{Post-training approximation and PWL-native training}}
\label{subsubsec:pwl_native}

The piecewise-linear activation defined above supports two distinct training
and deployment paths. In post-training approximation, the network weights
and biases are learned using smooth tanh activations. These activations are
subsequently replaced by piecewise-linear functions for optimization
embedding. Consequently, the embedded predictor can differ from the smooth
predictor whose error was minimized during training.

In PWL-native training, the piecewise-linear activation is used directly in
the training forward pass. The resulting weights and biases are therefore
fitted to the same activation representation used during prediction and
optimization. When the weights, biases, scaling, breakpoints, activation
values, and treatment of pre-activation domains are preserved, the
optimization embedding represents the trained predictor exactly, up to
numerical tolerances. No additional activation-replacement error is
introduced during embedding.

This distinction does not eliminate surrogate prediction error relative to
the underlying process. It eliminates the additional discrepancy caused by
changing activation functions after training. For PWL-native networks,
breakpoint selection determines the activation representation available
during fitting and the size of the subsequent optimization formulation; it
is not merely a post-training approximation setting.

PWL-native describes the activation representation rather than a particular
training algorithm. Gradient-based PWL-native training does not provide a
global training guarantee and should be distinguished from the restricted
SOSX-MAE training formulation introduced in the following section.
Similarly, using SOSX to optimize process inputs through a trained
PWL-native network does not imply that its weights were globally optimized.

\subsubsection{Mean Squared Error Formulations for Control Applications}
\label{sec:PWL-mse}
A major application of optimization-aware neural surrogates is model predictive control and dynamic optimization. In these applications, the optimization objective frequently involves tracking a desired process trajectory or minimizing deviations from a target operating condition. For such problems, the choice of loss function influences not only surrogate training but also the behavior of the resulting controller.

While mean absolute error (MAE) provides a linear objective that is robust to individual prediction errors, mean squared error (MSE) penalizes large deviations more heavily through its quadratic structure. Consequently, MSE is the more common objective in trajectory tracking and setpoint-control applications, where reducing large prediction errors is generally more important than minimizing the average absolute deviation. Controllers designed using MSE therefore tend to prioritize eliminating significant tracking errors, whereas MAE-based formulations typically produce more conservative behavior by treating all deviations uniformly.

Introducing an MSE objective changes the optimization formulation from a mixed-integer linear program to one containing a convex quadratic objective. This problem can be represented directly as a Mixed-Integer Quadratic Program (MIQP) or, equivalently, reformulated as a Mixed-Integer Second-Order Cone Program (MISOCP). Although these formulations describe the same convex optimization problem, their computational characteristics differ. MIQP solvers exploit quadratic programming algorithms, whereas conic solvers such as MOSEK exploit the underlying second-order cone structure through specialized interior-point methods, which can provide improved numerical robustness and computational performance for many convex quadratic optimization problems.

To support this broad range of optimization applications, PWLNet includes both MIQP and MISOCP reformulations for piecewise-linear neural surrogates. Together with the existing MILP formulations, users can readily investigate the trade-offs between linear, quadratic, and conic optimization models under a unified interface. This enables direct comparison of prediction accuracy, optimization model size, computational performance, and closed-loop control behavior while requiring only minimal modifications to the training and optimization workflow.

\subsection{SOSX-Based Piecewise-Linear Neural Network Training}
\label{subsec:sosx_training}

For the SOSX-based training path, PWLNet considers a single-hidden-layer neural network with a hyperbolic tangent activation in the hidden layer and a linear output layer. For a data point $(x_i,y_i)$, the network prediction is written as
\begin{equation}
\hat{y}_i
=
W_2 f(W_1x_i+b_1)+b_2,
\end{equation}
where $W_1$ and $b_1$ are the input-to-hidden weights and hidden-layer biases, $W_2$ is the hidden-to-output readout matrix, $b_2$ is the output bias, and $f(\cdot)=\tanh(\cdot)$.

Following the SOSX-based neural-network training formulation of K\"oksal, T\"urkay, and Ayd{\i}n, the output-layer weight matrix is not optimized simultaneously with the hidden-layer parameters \citep{koksal2026sosx_training}. Instead, the readout matrix is estimated from a preliminary conventional neural-network training step and then kept fixed during the SOSX optimization. This fixed readout is denoted by $W_2^\ast$. This assumption is closely related to the extreme learning machine assumptions, in which part of the network is fixed or precomputed so that the remaining training problem becomes simpler \citep{huang2006}.

With the fixed readout matrix, the SOSX training problem is formulated as
\begin{equation}
\min_{W_1,b_1,b_2}
\frac{1}{N}
\sum_{i=1}^{N}
\left|
W_2^\ast f'(W_1x_i+b_1)+b_2-y_i
\right|,
\label{eq:sosx_training_mae}
\end{equation}
where $f'(\cdot)$ denotes the piecewise-linear approximation of the hyperbolic tangent activation function. The trainable variables are therefore the input-to-hidden weights $W_1$, the hidden biases $b_1$, and the output bias $b_2$, while $W_2^\ast$ remains fixed.

The piecewise-linear approximation is defined over the breakpoint vector
\begin{equation}
a_0 < a_1 < \cdots < a_P,
\end{equation}
where $P$ is the number of linear pieces. For each segment $k=1,\ldots,P$, define the segment length and function increment as
\begin{equation}
u_k = a_k-a_{k-1},
\qquad
g_k = f'(a_k)-f'(a_{k-1}).
\end{equation}
For each data point, the hidden pre-activation is restricted to the approximation domain:
\begin{equation}
a_0 \leq W_1x_i+b_1 \leq a_P.
\end{equation}
Using incremental variables $Y_{i,k}$, the hidden pre-activation and the corresponding piecewise-linear activation value are represented as
\begin{align}
W_1x_i+b_1
&=
a_0+\sum_{k=1}^{P}Y_{i,k},
\label{eq:sosx_preactivation_incremental}
\
f'(W_1x_i+b_1)\\
&=
f(a_0)
+
\sum_{k=1}^{P}
\frac{g_k}{u_k}Y_{i,k}.
\label{eq:sosx_activation_incremental}
\end{align}

The absolute-value objective is linearized by introducing residual variables $e_i$ and nonnegative deviation variables $d_i$:
\begin{equation}
\begin{aligned}
e_i
&=
W_2^\ast f'(W_1x_i+b_1)+b_2-y_i,
&& i=1,\ldots,N,
\\
d_i &\geq e_i,
&& i=1,\ldots,N,
\\
d_i &\geq -e_i,
&& i=1,\ldots,N,
\\
d_i &\geq 0,
&& i=1,\ldots,N.
\end{aligned}
\end{equation}
The resulting SOSX training formulation can therefore be written as
\begin{equation}
\begin{aligned}
\min
\quad
& \frac{1}{N}\sum_{i=1}^{N}d_i
\\
\text{s.t.}\quad
& e_i
=
W_2^\ast f'(W_1x_i+b_1)+b_2-y_i,
&& i=1,\ldots,N,
\\
& d_i \geq e_i,
&& i=1,\ldots,N,
\\
& d_i \geq -e_i,
&& i=1,\ldots,N,
\\
& a_0 \leq W_1x_i+b_1 \leq a_P,
&& i=1,\ldots,N,
\\
& W_1x_i+b_1
=
a_0+\sum_{k=1}^{P}Y_{i,k},
&& i=1,\ldots,N,
\\
& f'(W_1x_i+b_1)
=
f(a_0)
+
\sum_{k=1}^{P}
\frac{g_k}{u_k}Y_{i,k},
&& i=1,\ldots,N.
\end{aligned}
\end{equation}

\begin{figure}
	\centering
		\includegraphics[width=1\textwidth]{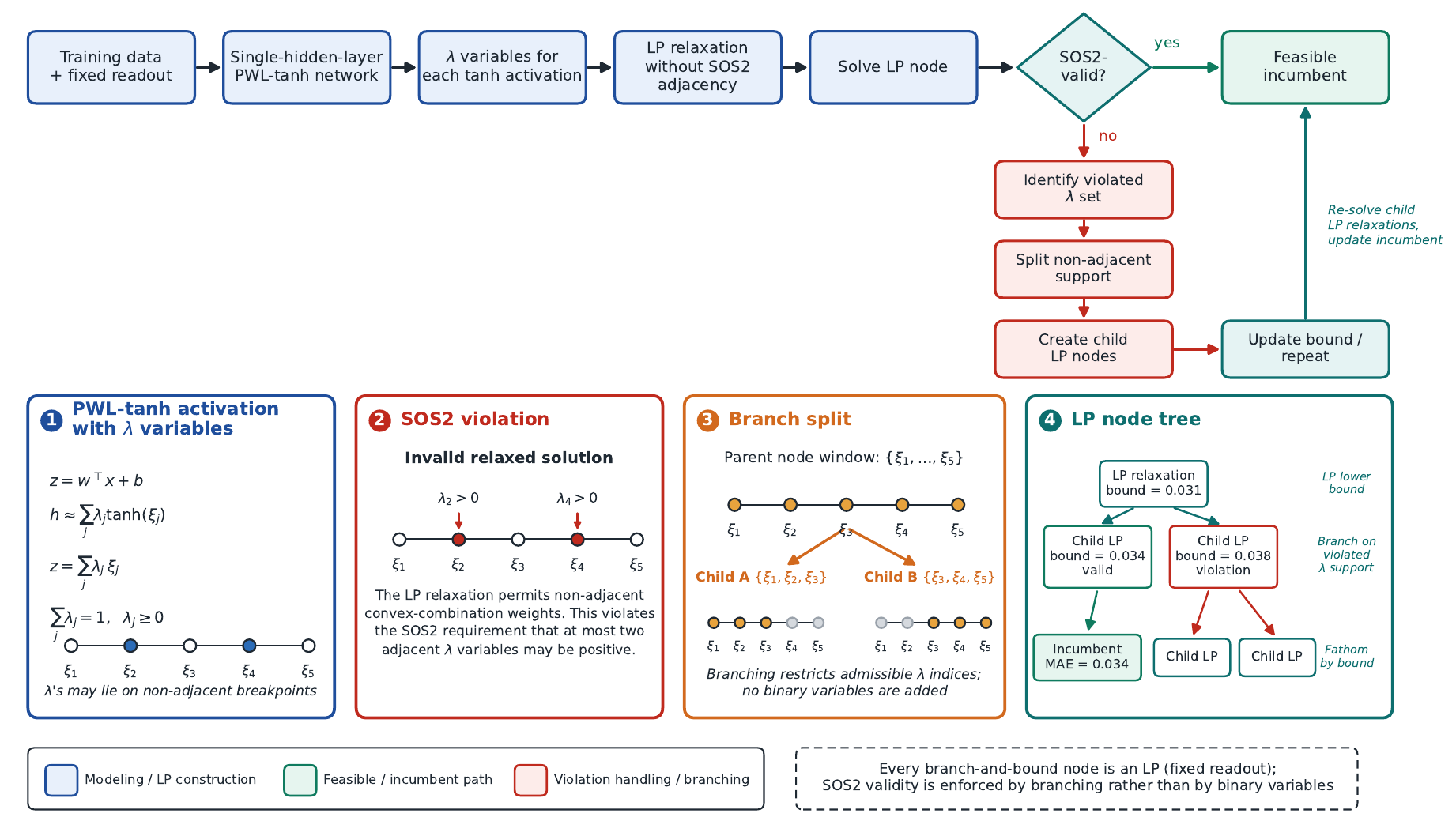}
	\caption{Schematic of the SOSX training for a single-hidden-layer PWL-tanh neural network. Each tanh activation is represented by convex-combination variables over fixed breakpoints. The relaxed training problem drops the SOS2 adjacency restriction and is solved as an LP because the output readout is fixed. If the LP solution activates non-adjacent $\lambda$ variables, the solution violates the SOS2 structure and therefore the brancher identifies the violated $\lambda$ auxiliary variable and splits the admissible breakpoint window into child nodes, and re-solves the corresponsing LP relaxations. SOS2-valid solutions provide feasible incumbent networks, while LP bounds are used to fathom subtrees. The procedure continues until the optimality gap, node limit, or time limit criterion is reached.}
	\label{fig:SOXtrain}
\end{figure}

The SOSX formulation enforces the valid piecewise-linear segment
structure through a customized branching procedure in the space of the
incremental variables $Y_{i,k}$, as illustrated in
Figure~\ref{fig:SOXtrain}. At each node of the branching tree, the additional node restrictions are linear equalities or inequalities on the $Y$ variables. Since $W_2^\ast$ is fixed, the prediction term
\begin{equation}
W_2^\ast f'(W_1x_i+b_1)+b_2
\end{equation}
is linear in the remaining decision variables after the tanh activation is replaced by its piecewise-linear approximation. The absolute-value loss is also represented using linear epigraph constraints. Consequently, each node subproblem in the SOSX branching algorithm is a linear program. The fixed-readout assumption is essential for preserving this LP structure. If the output weights were optimized together with the hidden-layer parameters, the prediction term
\begin{equation}
W_2 f'(W_1x_i+b_1)
\end{equation}
would contain products between trainable output weights and trainable hidden activations, introducing bilinear terms. Fixing $W_2^\ast$ removes this bilinear coupling and allows the SOSX branch-and-bound procedure to solve LP relaxations at each node rather than nonlinear or mixed-integer nonlinear subproblems.

\subsubsection{Incumbency Formulation}
\label{sec:sosx-incumbency}

The SOSX branch-and-bound search maintains an \emph{incumbent}, the best
SOS2-feasible solution found so far. A node solution is accepted as a new
incumbent only if (i) the node LP is feasible, (ii) no SOS2 adjacency
violation is detected among its convex-combination variables, and (iii) its
objective strictly improves the current incumbent by more than the
optimality tolerance. Because every node LP is a relaxation of its subtree,
the incumbent also drives fathoming: any open node whose relaxed bound
cannot improve the incumbent is discarded, and the difference between the
best open bound and the incumbent objective provides the absolute and
relative optimality gaps used as termination criteria.

In repeated-solve settings such as model predictive control, consecutive
optimization problems differ only through updated measured states or a
shifted horizon, so the previous solution is usually an excellent candidate
incumbent for the next problem. PWLNet therefore accepts a previous SOSX
result as a warm start. The previous \emph{objective value is never
trusted}, because the model data have changed. Instead, the warm start is
\emph{revalidated}: for every SOS set $s$, the stored convex-combination
weights $\lambda^{s}$ define their active support
\begin{equation}
  \mathcal{S}_s \;=\; \{\, k \;:\; \lambda^{s}_{k} > \varepsilon \,\},
\end{equation}
which is accepted only when it contains at most two adjacent breakpoints; the corresponding tight breakpoint window $w_s = [\min \mathcal{S}_s,\; \max \mathcal{S}_s]$ is then imposed on the \emph{updated} model and the resulting restricted problem is solved as a single LP. (For the incremental $Y$-variable formulation used in SOSX training, the active segment $k^{*}_s$ is identified as the first segment that is not filled to its length, $Y^{s}_{k} < u_{k}$, and the window is fixed to that single segment.) If this restricted LP is feasible and free of SOS2 violations, its objective seeds the incumbent before the root node is expanded, so large parts of the new tree are fathomed immediately; if it is infeasible or violates adjacency, the warm start is silently discarded and the search starts cold. This mechanism preserves global validity of the branch-and-bound search while transferring information between the related problems of a receding- or shrinking-horizon sequence.

\subsubsection{Persistent Solver Integration}

\label{sec:sosx-persistent}

A SOSX search may solve hundreds or thousands of node LPs, and consecutive nodes differ only in which convex-combination variables are fixed to zero. Rebuilding and re-transferring the full Pyomo model at every node would dominate the runtime. PWLNet therefore integrates the branching loop with persistent solver interfaces (Pyomo APPSI): the node relaxation is loaded into the solver once, and applying a breakpoint window merely fixes the $\lambda$ (or $Y$) variables outside the window to zero and releases the variables inside it. The relaxation object records exactly which variables changed fixing state between consecutive nodes and pushes only those updates to the solver, so each node reduces to an incremental LP re-solve that can exploit the solver's warm-started (dual simplex) basis.

Persistent mode is requested through \texttt{solver\_options=\{"persistent": True\}} and currently supports HiGHS (\texttt{appsi\_highs}) and Gurobi (\texttt{appsi\_gurobi}); common solver aliases are normalized to the corresponding APPSI interface. An infeasible node is reported back to the brancher as a fathomed node rather than raised as an exception, and solution loading is deferred until a feasible termination status is confirmed, which keeps the search loop robust when many restricted nodes are infeasible by construction. The non-persistent path remains available for any Pyomo-compatible LP/MILP solver, so persistence is an optimization, not a requirement.


\subsection{Optimization-Aware Simplification and Physics-Informed Training}
\label{sec:pruning-physics}

The size of the embedded optimization model depends directly on the size of the
neural network. Each hidden neuron adds variables and constraints to the
reformulation. PWLNet treats network compression and physics-informed training as
two complementary tools for producing surrogates. With these tools it is aimed that the formulation size is reduced and the interpretability of the model is increased without compromising too much from accuracy.

\subsubsection{Optimization-Aware Pruning}
\label{sec:pruning}

PWLNet applies pruning as a preprocessing step before reformulation, following
the workflow
\begin{equation}
\text{train network} \rightarrow \text{prune} \rightarrow \text{fine-tune}
\rightarrow \text{reformulate} \rightarrow \text{optimize}.
\label{eq:pruning-workflow}
\end{equation}

The iterative procedure formalized by \citet{lastrucci2026pruning} provides the methodological basis for this sequence. A dense network is first trained to convergence; its parameters or structural components are then assigned importance scores, low-scoring entries are removed by updating a pruning mask, and the score--prune cycle may be repeated until the target sparsity is reached. An optional fine-tuning stage subsequently adjusts only the retained parameters to recover predictive accuracy. PWLNet adopts this general separation between training, structural modification, and recovery, but uses it as a configurable preprocessing interface rather than prescribing a single iterative magnitude criterion. It therefore combines data-driven weight and neuron scores with activation-specific tests, such as dead or always-active ReLU detection and saturated-tanh detection, and evaluates the compressed network by prediction error, embedded formulation size, solver time, and true-process decision quality. The pruning methods considered in PWLNet are summarized in Table~\ref{tab:pruning-methods}.

\begin{table}[h]
\centering
\caption{Pruning methods considered in PWLNet.}
\label{tab:pruning-methods}
\begin{tabular}{lll}
\hline
Method & Activation & Purpose \\
\hline
Weight magnitude pruning & ReLU/tanh & Remove small weights. \\
Neuron magnitude pruning & ReLU/tanh & Remove weak hidden neurons. \\
Dead ReLU detection & ReLU & Remove neurons that are always inactive. \\
Always-active ReLU simplification & ReLU & Replace active ReLU units with linear equalities. \\
Saturated tanh detection & tanh & Identify tanh units operating near saturation. \\
Optimal-Brain-Damage-inspired pruning & ReLU/tanh & Remove low-saliency weights or neurons. \\
Group-lasso-style compression & ReLU/tanh & Drive whole neurons toward zero during training. \\
Fine-tuning after pruning & ReLU/tanh & Recover prediction accuracy after compression. \\
\hline
\end{tabular}
\end{table}

Weight magnitude pruning zeroes a weight $w_{ij}$ whenever
\begin{equation}
|w_{ij}| < \tau_w,
\label{eq:weight-magnitude}
\end{equation}
for a user-specified threshold $\tau_w$. Neuron magnitude pruning extends this
criterion to an entire neuron $j$ by comparing an importance score, such as
the norm of its incoming or outgoing weight vector $\mathbf{w}_j$, against a
neuron-level threshold $\tau_n$,
\begin{equation}
\lVert \mathbf{w}_j \rVert < \tau_n \;\Rightarrow\; \text{remove neuron } j.
\label{eq:neuron-magnitude}
\end{equation}

Dead ReLU detection and always-active ReLU simplification use the propagated
pre-activation bounds $L_j \le s_j \le U_j$. A ReLU neuron is dead and can be
removed if
\begin{equation}
U_j \le 0,
\label{eq:dead-relu}
\end{equation}
since $z_j = \max(0, s_j) = 0$ for every admissible input. Conversely, a
neuron is always active if
\begin{equation}
L_j \ge 0,
\label{eq:always-active-relu}
\end{equation}
in which case $z_j = s_j$ holds identically and the neuron can be represented
by a linear equality without an accompanying binary variable, removing one
source of combinatorial complexity from the embedded formulation.

Saturated tanh detection applies the same bound-propagation logic to tanh
units: if the propagated pre-activation range lies entirely outside the
piecewise-linear approximation domain, $[L_j, U_j] \not\subset (-a_P, a_P)$
for the outer breakpoint $a_P$, the unit is replaced by the constant
$z_j \approx \pm 1$, removing its nonlinear or piecewise-linear representation
from the optimization model.

Optimal-Brain-Damage-inspired pruning ranks weights or neurons by a
second-order saliency score using the diagonal Hessian approximation
\begin{equation}
\eta_i = \tfrac{1}{2} H_{ii}\, w_i^2,
\label{eq:obd-saliency}
\end{equation}
where $H_{ii}$ is the diagonal entry of the loss Hessian with respect to
parameter $w_i$; parameters with the smallest saliency are removed first,
limiting the resulting increase in training error. Group-lasso-style
compression instead acts during training or fine-tuning, adding a structured
sparsity penalty
\begin{equation}
\mathcal{R}_{\text{group}}(W) = \lambda_g \sum_{j} \lVert \mathbf{W}_{:,j} \rVert_2
\label{eq:group-lasso}
\end{equation}
to the training objective, where $\mathbf{W}_{:,j}$ collects the weights
associated with neuron $j$; driving entire columns of $W$ toward zero
identifies whole neurons for removal rather than individual weights.

After any of the above operations, a fine-tuning step retrains the remaining
parameters on the original supervised objective to recover the prediction
accuracy lost during compression. Each stage returns a \texttt{PruneOutcome}
record, and reporting utilities: compare parameter counts, prediction error,
and formulation size before and after compression.

The main hypothesis is that modest pruning can substantially improve
optimization tractability while preserving enough prediction accuracy for
downstream process optimization; this hypothesis is tested empirically in the
Williams--Otto and Van de Vusse case studies of Section~\ref{sec:benchmarks}.

\subsubsection{Physics-Informed Training}
\label{sec:physics-informed}

While network compression reduces the size of an already-trained surrogate,
physics-informed training changes what the surrogate learns in the
first place. In process systems engineering, partial physical knowledge, such
as mass balances, conservation relationships, operating bounds, or
monotonicity of specific outputs, is often available even when a complete
first-principles model is not. PWLNet incorporates this knowledge directly
into the training objective through the \texttt{PhysicsInformedTrainer}, so
that predictions violating known process behavior are penalized during
training rather than discovered only after optimization.

Given a supervised data loss $\mathcal{L}_{\text{data}}(W,b)$, such as the
mean squared or mean absolute error objectives used elsewhere in this section,
the physics-informed training problem augments the objective with a weighted
sum of differentiable residual terms,
\begin{equation}
\min_{W,b} \;\; \mathcal{L}_{\text{data}}(W,b)
+ \sum_{r=1}^{R} \lambda_r \, \mathcal{R}_r\bigl(x_i, \hat y_i; W, b\bigr),
\label{eq:physics-informed-loss}
\end{equation}
where each residual $\mathcal{R}_r$ measures the violation of a known process
relationship by the current prediction, and $\lambda_r \ge 0$ is a
user-specified penalty weight. Residual terms are supplied either through a
\texttt{PhysicsLoss} implementation or a declarative \texttt{PhysicsSpec},
which keeps the generic training loop independent of the process equations
used in a specific case study.

Two residual forms are used throughout the benchmark suite in Section~\ref{sec:benchmarks}.
Closure and conservation residuals penalize the extent to which predicted
outputs fail to satisfy a known balance relationship, such as a mass-fraction
closure constraint or a discretized ODE residual between consecutive
predicted states. Bound and monotonicity residuals penalize predictions that
fall outside known physical limits or violate a known qualitative trend, such
as the monotonic accumulation of an irreversible product over time. Both
residual forms are differentiable, so they can be minimized jointly with the
data loss using the same gradient-based training procedure used for the
unconstrained network.

Because physics residuals are evaluated only during training, they do not
change the mathematical form of the trained network, and therefore do not
affect the reformulation or embedding procedures described earlier in this
section. Physics-informed training is consequently compatible with every
formulation supported by PWLNet: it changes the parameters $W$ and $b$ that
are embedded, not the embedding itself \citep{koksal2023physics_pwl}.

\section{Results \& Discussions}
\label{sec:benchmarks}

The benchmark suite contains four cases with increasing relevance to process systems engineering which can be seen in Table~\ref{tab:cases}. The purpose is to test whether conclusions about linearized neural surrogates are consistent across generic regression, simple nonlinear process optimization, dynamic reactor behavior, and realistic process optimization.

\begin{table}[htbp]
\centering
\caption{Benchmark cases used in the study.}
\label{tab:cases}
\begin{tabular}{p{0.17\textwidth} p{0.34\textwidth} p{0.39\textwidth}}
\toprule
\textbf{Case} & \textbf{Role} & \textbf{Main evaluation focus} \\
\midrule
Boston housing & Small generic regression baseline & Prediction accuracy and formulation overhead. \\
Williams--Otto & Simple nonlinear process optimization & Surrogate embedding and objective quality. \\
Van de Vusse CSTR & Dynamic nonlinear reactor benchmark & Dynamic behavior and reactor optimization. \\
Hydroformylation & Realistic process optimization benchmark & Process constraints and economic optimization. \\
\bottomrule
\end{tabular}
\end{table}


\subsection{Case Study I: Regression Benchmark}
\label{sec:reg}

The Boston housing dataset (506 samples, 13 neighbourhood-level features,
one regression target \texttt{medv}) is used as a small generic regression
baseline and as the tutorial case for the complete PWLNet workflow:
training, piecewise-linear conversion, SOSX training, Pyomo export, and
embedded optimization. A 70/15/15 train/validation/test split is used, and
observed per-feature minima and maxima are pinned to every trained network
as input bounds so that the embeddings can derive tight pre-activation
intervals. Because several features are heavy-tailed, min--max input
scaling is used; this keeps the warm-start pre-activations inside the
bounded PWL-tanh domain $[-4, 4]$ required by the SOSX training
formulation.

\begin{figure}
\centering
\includegraphics[width=\textwidth]{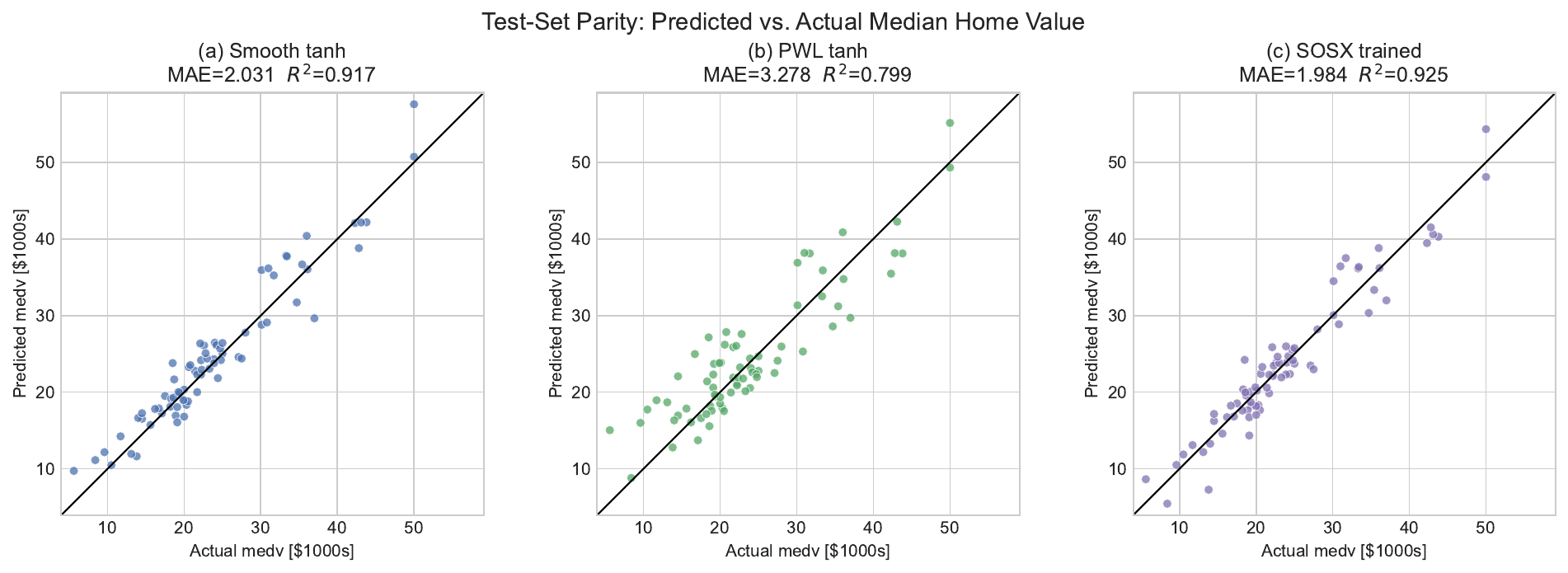}
\caption{Boston housing test-set parity plots for the smoothly trained tanh,
post-training PWL-tanh, and SOSX-trained PWL-tanh predictors sharing one
10-neuron single-hidden-layer architecture (notebook 13).}
\label{fig:boston-parity}
\end{figure}

Three pipelines share one architecture (13 inputs, 10 hidden tanh neurons,
1 linear output, 151 parameters). The smooth tanh network is trained with
Keras MSE; the post-training PWL-tanh predictor reuses the smooth weights and only replaces
the activation with the 5-breakpoint piecewise-linear approximation at
evaluation and embedding time; the SOSX-trained PWL-tanh network globally
re-optimizes the input layer and biases under the PWL-MAE objective with
the readout fixed, warm-started from the smooth model. Test-set results are summarized in Table~\ref{tab:boston-regression},
and the corresponding parity plots are shown in
Figure~\ref{fig:boston-parity}. The SOSX refinement
adds only 5.9\,s on top of the 69.2\,s warm-start fit and terminates at
global optimality after two LP nodes, slightly improving on the smooth
model in MAE terms (its native objective). All three
surrogates confirm that a compact single-hidden-layer tanh network is
sufficient for this benchmark (test $R^2$ up to $0.925$).

\begin{table}[htbp]
\caption{Boston housing test-set regression metrics and training cost for
the three pipelines (notebook 13). Metrics are evaluated using each row's predictor; the post-training PWL-tanh row uses the piecewise-linear
activation with the smoothly trained weights.}
\label{tab:boston-regression}
\centering
\begin{tabular}{lcccc}
\hline
Pipeline & MAE & RMSE & $R^2$ & fit time [s] \\
\hline
Smooth tanh predictor (Keras MSE) & 2.031 & 2.631 & 0.917 & 69.2 \\
Post-training PWL-tanh        & 3.278 & 4.105 & 0.799 & +0.0 \\
SOSX-trained PWL-tanh        & 1.984 & 2.511 & 0.925 & +5.9 \\
\hline
\end{tabular}
\end{table}

For the optimization track, each surrogate is embedded into a Pyomo model
whose 13 input features are decision variables bounded by the observed
data ranges, and the predicted median home value is maximized subject to
five planning-motivated constraints (air quality $\mathrm{nox} \le 0.50$,
crime $\mathrm{crim} \le 1.0$, school quality $\mathrm{ptratio} \le 18$,
$\mathrm{lstat} \le 15$, and average rooms $\mathrm{rm} \ge 6$). The post-training PWL-tanh / big-M MILP solves in 0.5\,s with HiGHS; the SOSX branch-and-LP path
solves the SOS2 embedding of the SOSX-trained PWL-tanh network in about 7\,s over
roughly 1{,}200 LP nodes. Both paths drive the interpretable levers in the
expected directions (maximum rooms, minimum \texttt{lstat}, minimum
pollution), and the binding constraints are exactly the planning limits.
A one-at-a-time sensitivity analysis around the constrained optimum ranks
\texttt{rm} and \texttt{lstat} as the dominant features, in agreement with
the classical regression literature, which supports the surrogate as a
plausible stand-in for the unknown ground-truth function. Since no
simulator exists for this dataset, this case validates usability and
consistency of the workflow rather than downstream decision quality.

\subsection{Case Study II: Williams-Otto Reactor}
\label{sec:wo}

The Williams--Otto benchmark models a continuous stirred-tank reactor with
reactions $A + B \to C$, $C + B \to P + E$, and $P + C \to G$. The two
operating decisions are the fresh-B feed rate $F_B \in [4, 7]$\,kg/s and
the reactor temperature $T_R \in [70, 100]$\,$^\circ$C; the surrogate maps
these to the six effluent mass fractions \citep{williamsotto,schmid2020williamsotto}. The optimization maximizes the
economic profit subject to the operating limits $X_A \le 0.12$ and
$X_G \le 0.08$; the known exact optimum is $F_B = 4.389$,
$T_R = 80.50$ with profit $75.82$\,\$/s, sitting exactly on both limits.
All surrogates are trained on 150 Latin-hypercube simulator samples
(60/20/20 split). Because the optimum is constraint-bound, the two
constraint outputs $X_A$ and $X_G$ are up-weighted ($8\times$--$10\times$,
renormalized to mean weight one) in the training loss; this consistently
places the learned feasibility boundary within a few $10^{-4}$ of the
exact one and is used in every Williams--Otto experiment below.

\begin{figure} 
\centering
\includegraphics[width=\textwidth]{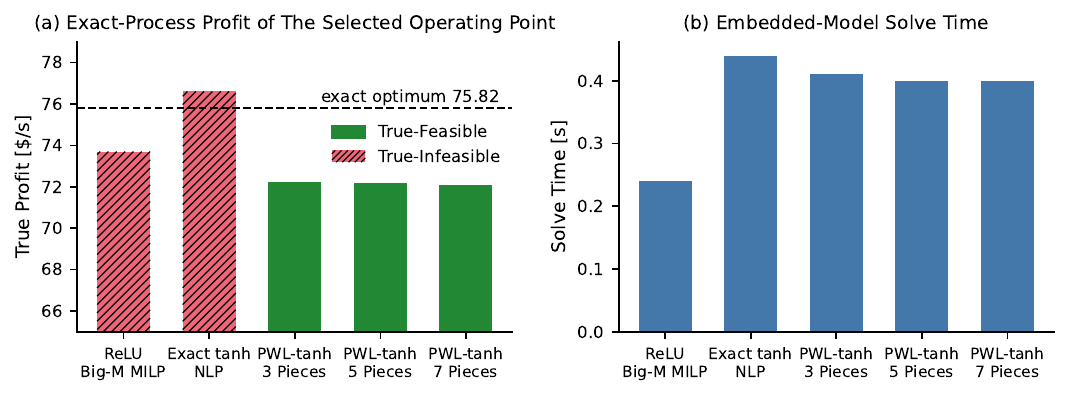}
\caption{Williams--Otto formulation comparison (notebooks 01--03).
(a) True-process profit at the operating point selected by each embedded
surrogate formulation; hatched bars violate a true operating limit
(the exact-tanh NLP exceeds the benchmark optimum only because its
selected point slightly violates the true $X_G$ limit). (b) Solve time of
the embedded model (HiGHS for the MILPs, Ipopt for the NLP).}
\label{fig:wo-formulations}
\end{figure}

\paragraph{Formulation comparison.}
Figure~\ref{fig:wo-formulations} compares the embedded formulations on
their native deployment paths. The exact tanh NLP (OMLT, Ipopt) is the
most accurate embedding: it selects an operating point essentially on top
of the exact optimum, but because the surrogate constraint boundary is
only accurate to a few $10^{-4}$, the selected point violates the true
$X_G$ limit by $\sim 5\times10^{-4}$; it is also only a local
solve. The ReLU big-M MILP (OMLT, HiGHS; 14 uncertain neurons, 14
binaries) reaches a true-process profit of $73.70$\,\$/s. The manual
post-training PWL-tanh / big-M MILPs deliver certified \emph{global} optima of the
embedded model in about $0.4$\,s for 3, 5, and 7 linear pieces, and all
three selected points are feasible for the true process, with true profits
of $72.1$--$72.2$\,\$/s. Two layers of approximation separate the
embedded optimum from the true process optimum: surrogate model error and
the PWL activation discretization; on this problem the remaining gap is
dominated by the economics relaxation and the surrogate error near the
constraint corner rather than by the breakpoint count. The bilinear
revenue terms $F_B X_P$ and $F_B X_E$ are handled by a piecewise McCormick
envelope (10 $F_B$ segments); a single-segment McCormick envelope
overstates the achievable profit by roughly $30$\,\$/s and misplaces the
recommended operating point, so economics relaxation quality matters as
much as the network formulation. Finally, the same PWL block supports MAE
(MILP), quadratic (MIQP), and second-order-cone (MISOCP) objectives: on an
inverse-design problem matching a target composition, the MIQP and MISOCP
variants recommend operating points within $0.01$ of each other (solved by
MOSEK in under $0.1$\,s), confirming that the problem class is set purely
by the objective formulation.

\begin{figure}
\centering
\includegraphics[width=\textwidth]{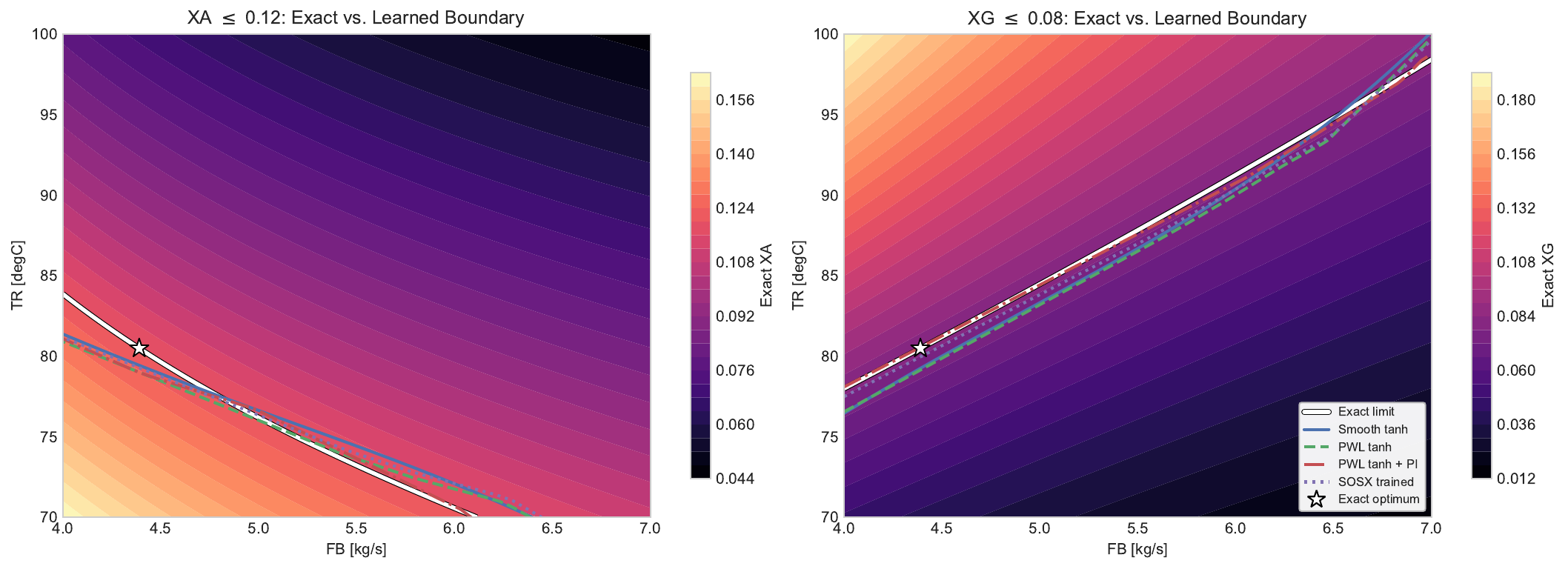}
\caption{Learned versus exact Williams--Otto operating constraints
(notebook 08). Filled contours show the exact process outputs; solid
white curves are the true limits $X_A = 0.12$ and $X_G = 0.08$, and
the dashed curves are the boundaries implied by each surrogate's
predicted outputs. The embedded optimizer enforces the limits on the
surrogate outputs, so the offset between dashed and solid curves is
the feasibility error that every embedded formulation inherits,
independently of the solver.}
\label{fig:wo-constraints}
\end{figure}

The feasibility boundaries that each embedded optimizer actually
enforces are compared with the exact ones in
Figure~\ref{fig:wo-constraints}; the gap between the learned and exact
$X_G$ boundary is what converts a surrogate-feasible optimum into a
marginally infeasible true operating point.

\begin{figure} 
\centering
\includegraphics[width=\textwidth]{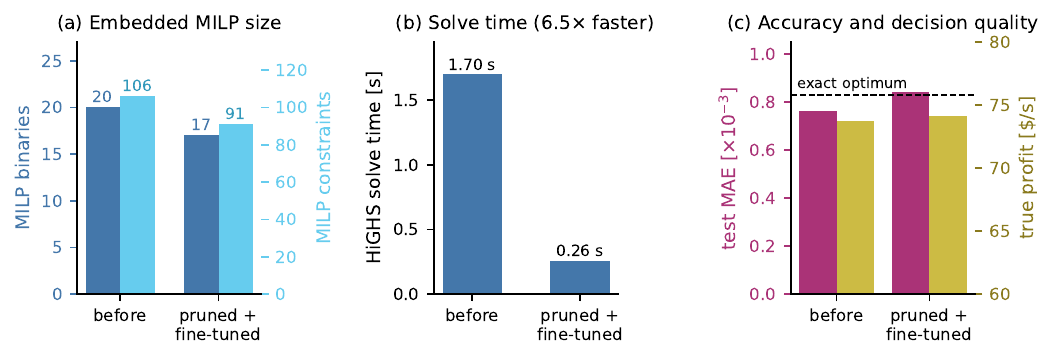}
\caption{Optimization-aware ReLU pruning on Williams--Otto (notebook 05).
(a) Embedded MILP size before and after pruning and fine-tuning;
(b) measured HiGHS solve time; (c) test MAE and true-process profit of
the recommended operating point.}
\label{fig:wo-pruning}
\end{figure}

\paragraph{Optimization-aware pruning.}
Figure~\ref{fig:wo-pruning} shows the effect of the compression pipeline
on a $[10,10]$ ReLU surrogate embedded as a big-M MILP. Exact bound-based
simplification removes nothing on this well-trained network; gentle
magnitude pruning (10\% of weights) followed by neuron-importance pruning
(10\% of neurons) and constraint-weighted fine-tuning removes 3 of 20
hidden neurons. Because every uncertain ReLU neuron contributes one binary
variable, this is direct binary-variable pruning: the embedded MILP shrinks
from 20 to 17 binaries and 106 to 91 constraints, and the measured HiGHS
solve time drops from $1.70$\,s to $0.26$\,s ($6.5\times$), while the test
MAE moves only from $7.6\times10^{-4}$ to $8.4\times10^{-4}$ and the
true-process profit of the recommended operating point is retained
($73.7 \to 74.1$\,\$/s). The upstream costs are modest as well (117\,s
initial training, 18\,s compression and fine-tuning), supporting the
central pruning hypothesis: the embedded model becomes smaller and faster
to solve while the decision it recommends stays essentially unchanged.

\begin{figure} 
\centering
\includegraphics[width=\textwidth]{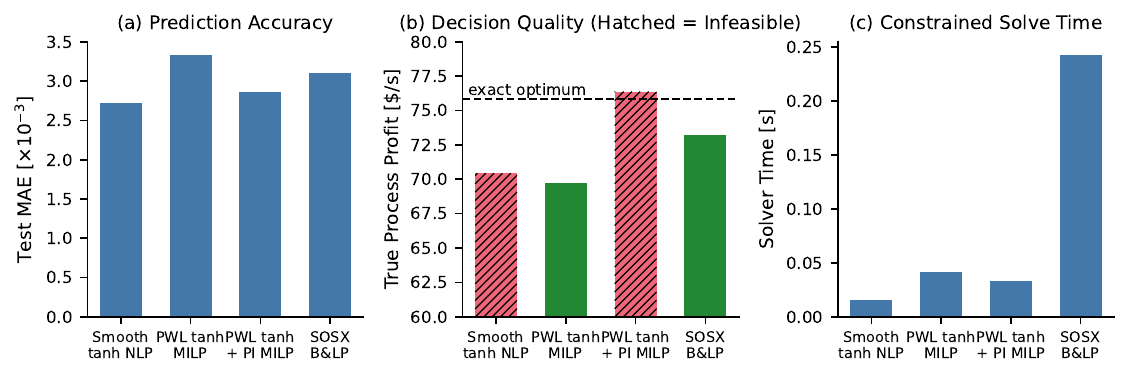}
\caption{Williams--Otto training-pipeline comparison 
(notebook 08): Exact tanh / IPOPT, post-training PWL-tanh / big-M,
physics-informed post-training PWL-tanh / big-M, and SOSX-trained
PWL-tanh / SOSX. (a) Test MAE for each deployed predictor; (b) true-process profit
at the selected operating point, hatched where the point violates a true
operating limit; (c) solver time of the constrained solve.}
\label{fig:wo-pipelines}
\end{figure}
\paragraph{Physics-informed training.}
Adding differentiable mass-fraction closure and operating-bound penalties
to the training loss of a $[16,16]$ tanh surrogate costs essentially
nothing in accuracy (test MAE $4.3\times10^{-4}$ plain versus
$4.6\times10^{-4}$ physics-informed) and reduces the mean normalized
steady-state balance residual by roughly 20\%
($4.3\times10^{-3} \to 3.5\times10^{-3}$). Under the exact-tanh NLP
embedding both models land on the true optimum (true profits 75.93 and
75.88\,\$/s versus the 75.82 benchmark, with only a
$\sim\!10^{-4}$ $X_G$ overshoot); under the 7-piece PWL/MILP embedding
both remain strictly feasible, and the physics-informed model is the more
conservative of the two. At this problem scale the downstream error is
dominated by the PWL discretization rather than by the training variant.

\paragraph{SOSX training and the four-pipeline comparison.}
Figure~\ref{fig:wo-pipelines} compares four end-to-end pipelines on one
fixed split with a single-hidden-layer, 6-neuron architecture: Exact tanh / IPOPT
(Keras MSE), post-training PWL-tanh / big-M (same smooth-trained weights),
physics-informed post-training PWL-tanh / big-M, and SOSX-trained
PWL-tanh / SOSX (global PWL-MAE training with fixed readout). Prediction
accuracy is comparable across all four deployed predictors (test MAE
$2.7$--$3.3\times10^{-3}$); the SOSX refinement itself costs only 1.2\,s
on top of its warm-start fit and terminates at global optimality after two
LP nodes. On the constrained economic optimization, the two formulations
that enforce the operating limits through the \emph{same} piecewise-linear
model that they optimize (post-training PWL-tanh / big-M and SOSX-trained PWL-tanh / SOSX) select truly feasible
operating points, whereas the smooth-NLP and PI variants land marginally
outside the true feasible region; the SOSX-trained PWL-tanh network attains the
best feasible true profit ($73.2$\,\$/s). Under training-input noise
(5\% of each input range, notebook 10), all pipelines degrade gracefully:
test MAE roughly doubles at worst, no formulation becomes unstable, and
the noisy-input SOSX solve still terminates with certified gaps, although
its selected operating point drifts to a truly infeasible corner --- a
reminder that constraint fidelity of the surrogate, not the optimizer,
is the binding safety factor.

\subsection{Case Study III: Van de Vusse Reactor}
\label{sec:VdV}

\subsubsection{Problem Formulation}
\label{sec:vdv-problem}

The Van de Vusse benchmark is the isothermal CSTR with reaction
network $A \to B \to C$ and $2A \to D$, states
$x = [C_A, C_B, C_C, C_D]$ and the dilution rate
$D = F/V \in [0.1, 4.0]$\,min$^{-1}$ as the manipulated input. The
control task is an operating-point transition with fixed final time:
starting from the loaded inverse-response steady state
$x(0) = [3.0,\, 1.117,\, 3.258,\, 1.3125]$\,mol/L at
$D_0 = 4/7$\,min$^{-1}$, maximize the terminal concentration of the
desired product $C_B$ at $t = 10$\,min using a piecewise-constant
input sequence $D_0, \dots, D_{99}$ with one move every 6\,s. The
reference solution is computed by coarse-to-fine direct shooting on
the true ODE and reaches a terminal $C_B$ of $1.614$\,mol/L --- well
above the best steady-state value of $1.266$\,mol/L --- through an
end-loaded strategy: drive $D$ high to accumulate reactant $A$, then
crash $D$ to its lower bound so the accumulated $A$ converts into a
terminal $C_B$ spike.

\begin{figure} 
\centering
\includegraphics[width=\textwidth]{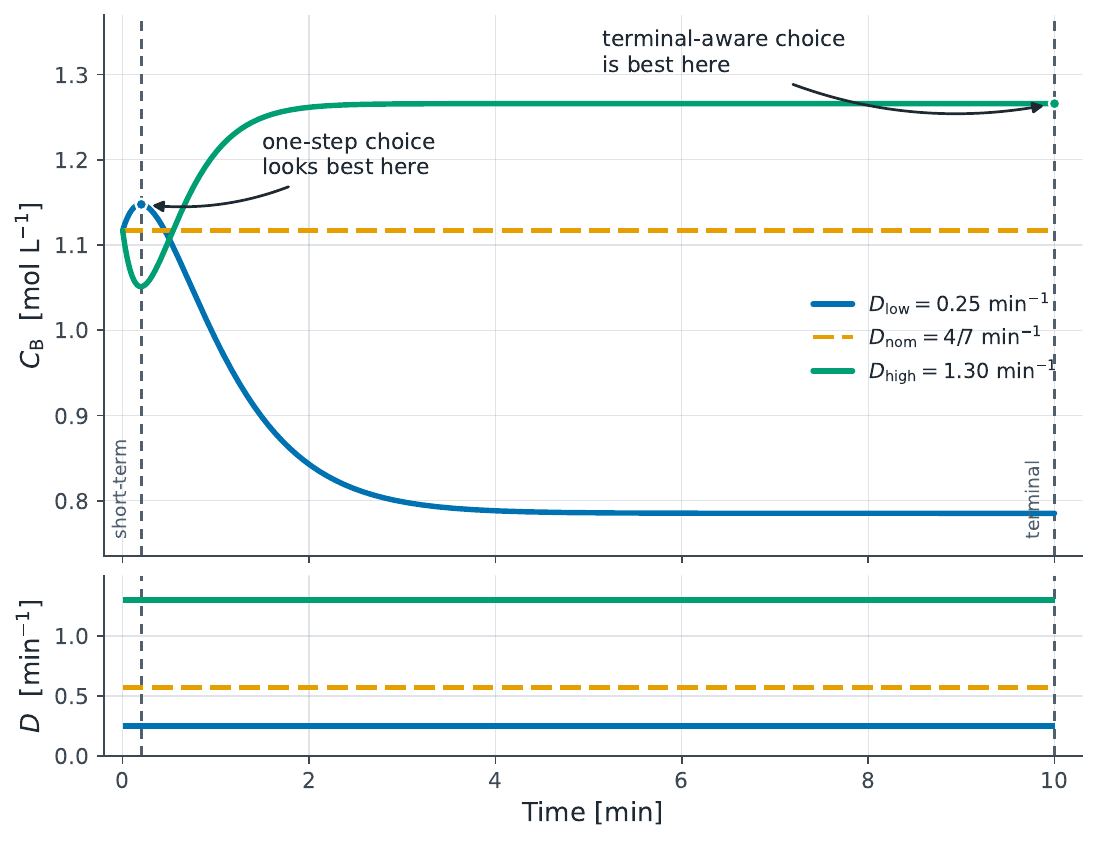}
\caption{Inverse-response motivation for terminal-aware Van de Vusse
control. Starting from the same nominal steady state at
$D_\mathrm{nom}=4/7$\,min$^{-1}$, three constant dilution rates are held
(bottom panel). The low rate raises $C_B$ first and looks best on a short
horizon, whereas the high rate initially depresses $C_B$ (inverse
response) but delivers the highest terminal value, so one-step-greedy and
terminal-aware controllers select opposite moves.}
\label{fig:vdv-motivation}
\end{figure}

Because the process exhibits inverse response, a one-step-greedy
controller that maximizes $C_{B,k+1}$ selects poor moves
(Figure~\ref{fig:vdv-motivation}). 
The learning problem is therefore formulated as \emph{terminal-aware,
shrinking-horizon output prediction}: the surrogate input is the
measured state, the next $P$ candidate moves, and a normalized horizon
feature,
\begin{equation}
  \big[C_{A,k},\, C_{B,k},\, C_{C,k},\, C_{D,k},\,
       D_k, \dots, D_{k+P-1},\, m/P\big]
  \;\longmapsto\;
  \Delta C_B^{(m)} = C_{B,k+m} - C_{B,k},
  \qquad m = \min(P,\, H - k),
\end{equation}
with maximum look-ahead $P = 10$ steps (1\,min). The active horizon
$m$ shrinks as the run approaches the fixed terminal time, so near the
end the surrogate predicts $C_B$ at the true terminal step rather than
a constant $P$ steps beyond it; the training windows cover every step
$k = 0, \dots, H-1$ so the terminal exit arc is represented in the
data.
The receding controller fixes the current measured concentrations as
equality constraints, maximizes the predicted end-of-horizon $C_B$
minus a small move penalty over the bounded future input sequence,
applies only the first move to the true plant, and repeats. Without
the shrinking horizon and terminal data coverage, a fixed-$P$
controller has no terminal-time awareness, keeps $D$ pinned high near
the end, and stalls at a terminal $C_B$ of roughly $0.99$ versus the
$1.61$ reference.

\begin{figure} 
\centering
\includegraphics[width=\textwidth]{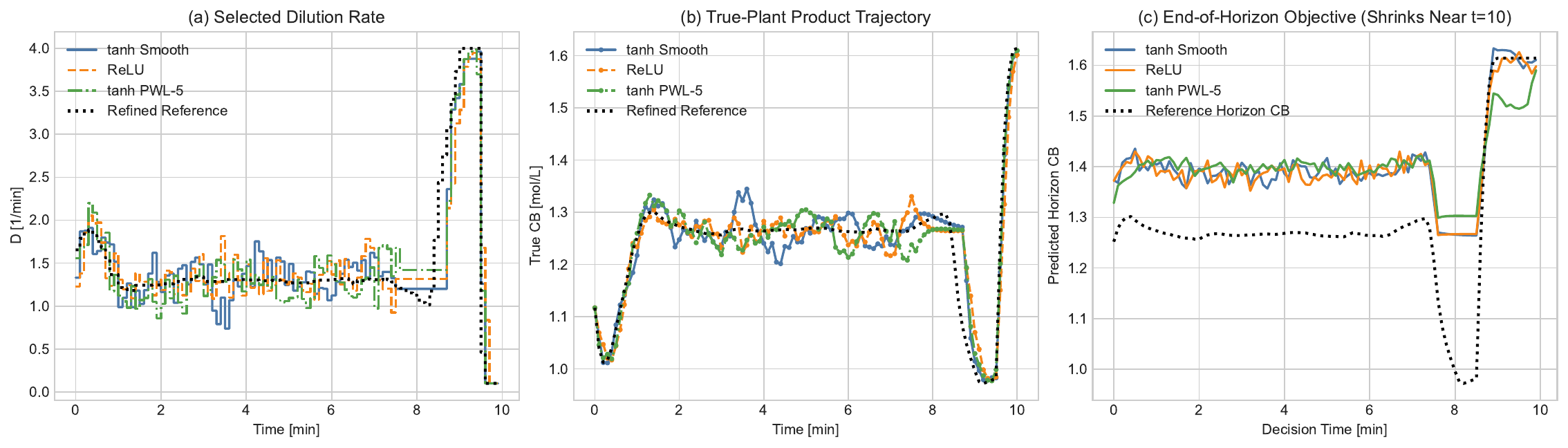}
\caption{Van de Vusse closed-loop receding-horizon control with the
trajectory-trained surrogates (notebook 14). (a) Dilution-rate sequences
selected by the Exact tanh / NLP, ReLU / big-$M$, and post-training
PWL-tanh / SOS2 controllers against the refined
direct-shooting reference; (b) true-plant $C_B$ trajectories; (c)
predicted end-of-horizon $C_B$, whose active horizon shrinks toward the
fixed terminal time.}
\label{fig:vdv-closedloop}
\end{figure}

Three horizon-aware surrogates are trained on 2{,}500 samples
generated from windows on the reference trajectory plus local state
and input perturbations: a smoothly trained tanh network, a ReLU network, and the
post-training PWL-tanh predictor using the same tanh weights (5
breakpoints). Test MAEs on the horizon target are $6.4\times10^{-3}$
for the smooth tanh predictor, $7.4\times10^{-3}$ for ReLU, and
$4.4\times10^{-2}$ for the deployed post-training PWL-tanh predictor. In
closed loop (Figure~\ref{fig:vdv-closedloop}), all three controllers execute the
build-up-then-crash exit arc and reach terminal $C_B$ values of
$1.6096$, $1.6017$, and $1.6099$\,mol/L respectively (99.2--99.7\% of the direct-shooting reference , $1.6142$) with all
state bounds respected and total closed-loop computation of a few
seconds for the complete 100-move run.

\subsubsection{Validation of the Learned Control Policy}
\label{sec:vdv-validation}

The surrogate-based controllers above closely match the
direct-shooting reference. However, this agreement alone does not
demonstrate that the surrogate has learned the underlying reactor
dynamics. An alternative explanation is that the optimization succeeds
simply because the surrogate was trained using data generated from the
optimal control trajectory itself, effectively allowing it to memorize
that particular solution. To distinguish between these possibilities,
two complementary validation studies are performed.

\paragraph{Global state-space training.}
In the first study, the surrogate is trained using samples distributed
throughout the reachable reactor state space rather than only along
the optimal trajectory. Since generating optimal trajectories from
many initial conditions would be computationally prohibitive, a
relatively coarse sampling strategy is employed. Although the
resulting dataset contains fewer samples than would typically be
required for high-resolution state-space coverage, it exposes the
surrogate to a much broader range of reactor states and control
actions. This experiment evaluates whether the surrogate learns the
fundamental reactor dynamics and can recover near-optimal control
policies without ever being trained on the optimal trajectory
corresponding to the optimization problem under consideration.

\paragraph{Local trajectory enrichment.}
The second study focuses on the original control problem with the same
fixed initial condition used throughout the benchmark. Instead of
constructing a global dataset, the training data are enriched locally
by augmenting the optimal trajectory with multiple suboptimal
trajectories generated from alternative control sequences. These
additional trajectories expose the surrogate to a larger portion of
the local state space surrounding the optimum while maintaining the
same optimization objective and initial condition. If the optimizer
consistently converges to the same solution despite being trained on a
richer set of trajectories, this indicates that the surrogate has
learned the local structure of the optimization landscape rather than
simply reproducing the single trajectory used for training. Together,
the two studies investigate surrogate generalization at different
scales: global reactor dynamics across the state space, and the local
optimization landscape of a specific operating scenario.

\paragraph{Results.}

\begin{figure} 
\centering
\includegraphics[width=\textwidth]{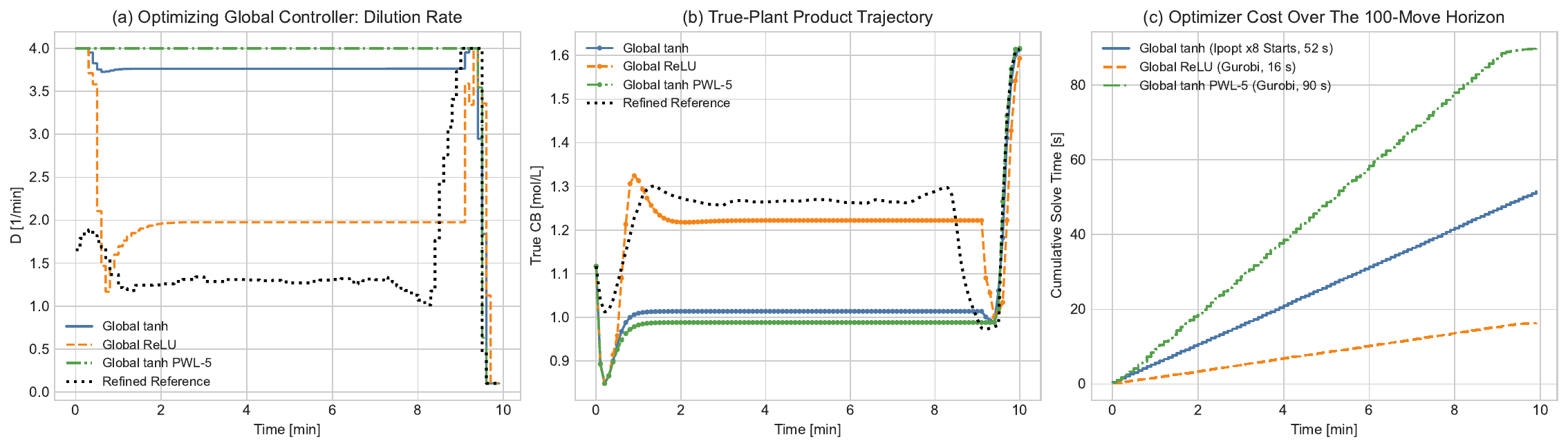}
\caption{Validation with globally trained surrogates and per-move
optimizing controllers (notebook 14). (a) Dilution profiles obtained by
optimizing directly through Exact tanh / IPOPT (multistart), ReLU /
big-$M$ (Gurobi), and post-training PWL-tanh / SOS2 (Gurobi) surrogates,
none of which were trained on the optimal trajectory; (b) true-plant $C_B$ trajectories against the reference; (c)
cumulative solver cost over the 100-move horizon.}
\label{fig:vdv-global_validation}
\end{figure}

In the global state-space study, the surrogate is retrained on 13{,}650
samples drawn uniformly from the full state--input box (never touching
the optimal trajectory), and at each step the receding controller
\emph{optimizes} the future dilution moves directly through the
surrogate rather than ranking a fixed candidate list: an Exact tanh / IPOPT NLP (multistart), a ReLU / big-$M$ MILP (Gurobi),
and a post-training PWL-tanh / SOS2 MILP (Gurobi). The tanh surrogate
trained on the global state-space samples predicts the horizon target along the
(never seen) optimal path with $R^2 = 0.996$, and all three optimizing
controllers recover the build-up-then-crash exit arc, reaching terminal
$C_B$ values of 100.0\% (tanh), 98.7\% (ReLU), and 100.1\% (PWL-5) of the
reference (Figure~\ref{fig:vdv-global_validation}a); the Exact tanh / IPOPT and post-training PWL-tanh / SOS2 controllers
marginally exceed the fixed direct-shooting reference. Since
neither the training data nor the optimizer is given the optimal input
sequence, this performance can only come from a learned representation
of the reactor dynamics. The price of breadth over depth is a lower
integral $C_B$: the tanh and PWL-5 controllers hold the reactor at a
leaner mid-horizon state (integral $C_B \approx 10.3$ and $10.1$ for
the Exact tanh / IPOPT and post-training PWL-tanh / SOS2 controllers,
respectively, versus $12.4$ for the reference and the trajectory-trained
controllers) before the terminal spike. The local enrichment study
confirms the complementary property: when the training windows mix the
optimal trajectory with many suboptimal perturbed trajectories, the
receding controller continues to converge to essentially the same
solution, indicating that the surrogate represents the local
optimization landscape rather than memorizing a single input sequence.

\subsubsection{Recurrent ReLU Surrogate MPC}
\label{sec:vdv-rnn}
The same reactor is used to demonstrate the recurrent engine of
Section~\ref{subsec:relu_rnn_mip}. 
A physics-informed ReLU RNN is trained on sequence data (control
profiles $D[t]$ mapped to state trajectories) with two penalty terms
alongside the data loss: a high-weight bound penalty keeping predicted
concentrations inside their physical limits, and a discretized Van de
Vusse ODE residual that forces the predicted trajectories to obey the
known mass balances. The cell is deliberately oversized (32 hidden
units) and then pruned: exact dead-unit removal (zero units on this
trained cell) followed by importance-based trimming to 12 units, which
increases the open-loop validation MAE on $C_B$ from $0.019$ to
$0.056$\,mol/L --- a modest cost for a 62\% reduction of the MILP
binaries per time step. Time-dependent interval bounds are then
propagated over the horizon and the pruned cell is unrolled into a
receding-horizon \texttt{relu\_bigm} MILP.

The unrolled embedding is tested for exactness. The RNN prediction values and the MILP formulation output is tested against each other. The result show that there is no formulation error and MILP formulation match RNN prediction to machine precision ($\sim10^{-16}$).

To the unrolled MILP model an economic objective is given for the final problem formulation. In this formulation, over a 20-step horizon (2 min) the solver must find the dilution profile that maximizes a running payoff.

\begin{equation*}
  J \;=\; \Delta t \sum_{t} \bigl( p_B\,C_{B,t} - c_C\,C_{C,t}
        - c_D\,C_{D,t} - c_F\,D_t \bigr)
        \;-\; w \sum_{t} \lvert \Delta D_t \rvert ,
\end{equation*}

Product revenue net of byproduct-disposal and feed costs is subjected to an actuator ramp-rate limit $\lvert \Delta D_t \rvert \le 0.25$ per move, so the optimizer must plan a ramp rather than jump to the best operating point. The formulation stays a pure MILP: the payoff is linear in the unrolled variables, and the $\lvert \Delta D_t \rvert$ auxiliaries both linearize the move penalty and, upper-bounded by the ramp limit, enforce the actuator constraint. 
A 20-step chain of ReLU blocks is, however, a hard MILP. The propagated interval boxes grow geometrically with the horizon, so the late-step big-M constants balloon and an untuned solve can run for hours. Two exact tightenings restore tractability: ramp-aware bound propagation, which threads the interval propagation one step at a time using only the dilution interval the actuator can actually reach (shrinking the terminal hidden box by an order of magnitude, $520$ vs. $6018$), and control-move blocking over pairs of steps (ten decision moves over twenty steps). Gurobi then certifies the optimum to a 0.9\% gap in about three and a half minutes (240 binaries), and the plan replayed through the NumPy cell matches the MILP outputs to solver feasibility tolerance ($\sim\!10^{-6}$).

The unrolled MILP earns $J = 3.40$ against $J = 3.27$ for the best ramp-to-and-hold policy, a $4.0\%$ economic gain from dynamic operation alone. The optimal plan ramps the dilution up to $D \approx 1.55$, glides back down, and ends with a terminal feed cut to the lower bound that harvests product concentration while saving feed cost; along the way the plant’s transiently reaches $C_B \approx\!1.37$\,mol/L, well above the steady-state ceiling of $1.266$ mol/L, exploiting the same inverse-response mechanism as the reference solution of Section~\ref{sec:vdv-problem}. Notably, the pruned surrogate under-predicts this transient peak ($\approx\!1.27$ vs. $1.37$\,mol/L on the plant), yet its plan still outperforms every ramp-to-and-hold policy on the true process. This confirms that compact recurrent ReLU surrogates are a viable exactly-embeddable alternative to the feedforward horizon-output surrogates when a full dynamic state model is desired.

\begin{figure} 
\centering
\includegraphics[width=\textwidth]{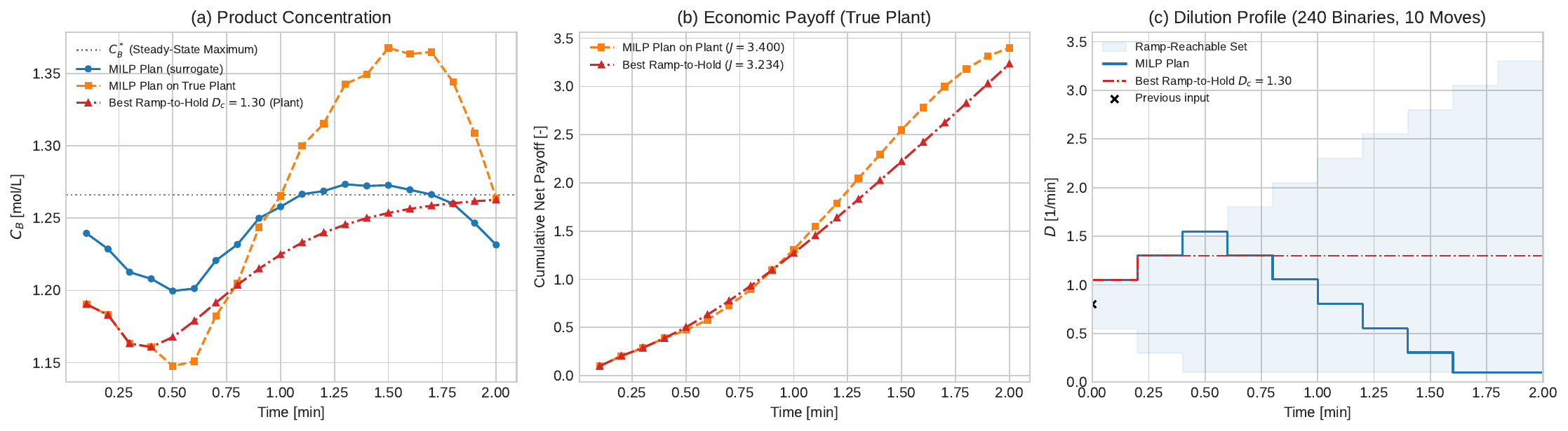}
\caption{Economic MPC with the pruned, physics-informed ReLU RNN
surrogate on the Van de Vusse reactor (notebook 12). Left: product
concentration under the optimal MILP plan on the surrogate and the true
plant against the best ramp-to-and-hold baseline; center: cumulative
economic payoff realized on the plant; right: blocked dilution profile
inside the ramp-reachable envelope.}
\label{fig:vdv_rnn}
\end{figure}

Figure~\ref{fig:vdv_rnn} summarizes the study: the product
concentration under the optimal plan on the surrogate and the true
plant against the baseline (left), the cumulative economic payoff on
the plant (center), and the blocked dilution profile inside the
ramp-reachable envelope (right).

\subsection{Case Study IV: Hydroformylation Reactor}
\label{sec:Hform}
\begin{figure}
\centering
\includegraphics[width=0.85\textwidth]{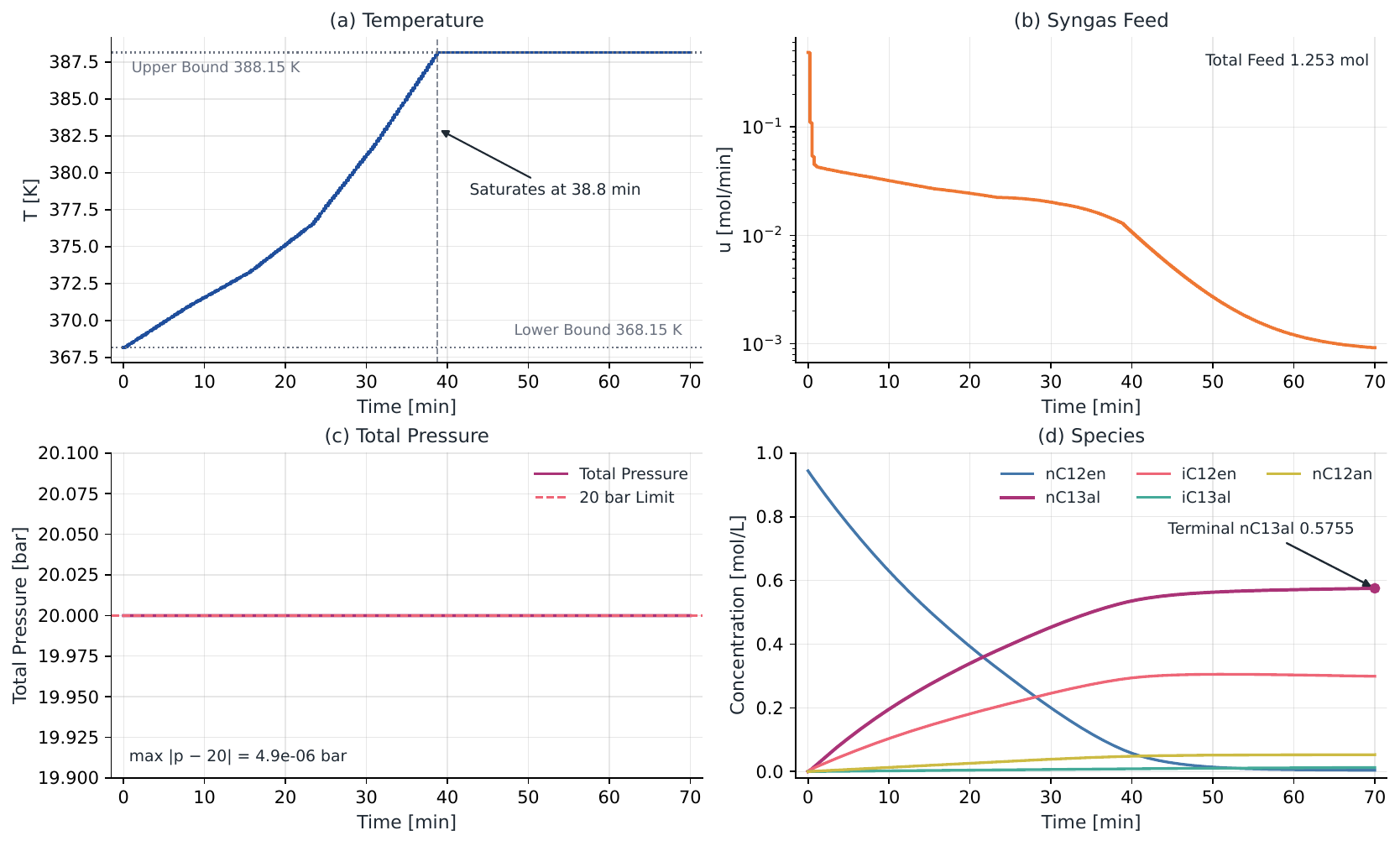}
\caption{Certified reference solution of the hydroformylation
optimal-control problem: multiple shooting on the true stiff ODE at
15-second resolution yields a pressure-holding policy that maintains the
total pressure at the 20-bar limit while ramping the temperature to its
upper bound and then holding it, reaching a terminal $n$C13al of $0.5755$
after the 70-minute batch.}
\label{fig:hydro-reference}
\end{figure}

The hydroformylation of 1-dodecene in a thermomorphic solvent system is
used as the realistic process optimization benchmark. The plant is a
stiff semi-batch reactor model with nine states: dissolved gas
concentrations $c_{H_2}$ and $c_{CO}$, five liquid-phase species including
the desired linear aldehyde $n$C13al, and the two gas partial pressures
$p_{H_2}$ and $p_{CO}$ which is based on the kinetics of Kiedorf et al.\ and
the reaction-route analysis of Hentschel et al.\
\shortcites{kiedorf2014kinetics,hentschel2015route}%
\citep{kiedorf2014kinetics,hentschel2015route}. The manipulated inputs are
the reactor temperature $T \in [368.15, 388.15]$\,K and the synthesis-gas
feed rate $u \in [0, 7]$\,mol/min. The embedded controllers restrict the
feed to $u \leq 0.4852$\,mol/min, the reference peak and the range supported
by their training data. The operating constraint is the total
pressure $p_{H_2} + p_{CO} \le 20$\,bar, and the objective is the terminal
amount of $n$C13al after the 70-minute batch. The reference optimum is a
multiple-shooting solution of the true ODE at 15-second resolution: a
pressure-holding policy that maintains the total pressure at the 20-bar
limit while ramping the temperature to its upper bound and then holding
it, reaching a terminal
$n$C13al of $0.5755$ (Figure~\ref{fig:hydro-reference}). The dynamic
optimization and uncertainty-aware NMPC setup follows the semi-batch
hydroformylation studies of \citet{aydin2017dynamic,aydin2018nmpc}.

The controller study builds on the transition-surrogate and
physics-informed experiments in notebooks 04 and 07, with the final
comparison taken from notebook 09. The full-state transition map
$[\tau,x_k,T_k,u_k]\mapsto x_{k+1}$ is trained on 1{,}230 simulated
trajectories, giving 344{,}400 transitions at 15-second resolution. The
30 baseline trajectories are supplemented by 600 trajectories with
randomized control timing and 600 with temperature profiles concentrated
near the input bounds, so that the training data cover more than small
perturbations of the reference policy. Plant parameters are sampled for
each trajectory, and transitions are split randomly into 80\% training
and 20\% test data. The results use the saved networks from seed 2479085,
selected in an earlier small-dataset sweep, and therefore represent a
single-seed comparison.

The alternatives include plain and physics-informed tanh and ReLU
networks, neuron pruning followed by fine-tuning, and PWL-native training.
The plain networks have 16 hidden neurons (361 parameters), while pruning
20\% of the neurons leaves 13 (295 parameters). Training uses 300 epochs,
with 50 additional epochs for fine-tuning and a physics-penalty scale of
0.2 for the physics-informed variants. The final PWL-native model uses
8 hidden neurons (185 parameters). Its three-piece activation, with
breakpoints $[-4,-1,1,4]$, is retained during training and SOSX embedding.
The smoothly trained tanh networks are evaluated both through the same
post-training PWL replacement with SOSX and through their exact smooth
activation with IPOPT. ReLU networks use an exact big-$M$ MILP embedding.

The development study varied the prediction horizon, activation
breakpoints, pressure margins, gas-inventory valuation, and training-data
coverage. Increasing the horizon did not consistently improve the realized
policy: some earlier PWL-native configurations benefited from $H=3$,
whereas the final compact model performs well with $H=1$. Longer prediction
chains also propagate surrogate error and, for post-training PWL-tanh,
activation-replacement error. Similarly, increasing the number of
breakpoints or imposing larger pressure margins did not consistently
improve the closed-loop result. Valuing the remaining gas inventory helped
avoid premature feed reduction, while broadening the training trajectories
improved coverage of the controls selected by the optimizer. These trials
motivate selecting the horizon together with the surrogate, rather than
treating a longer horizon as an improvement by itself.

In the final comparison, each short-horizon controller plans one 15-second
move ($H=1$), applies it to the true stiff-ODE plant, and repeats this for
280 moves. All start at $368.15$\,K and enforce a temperature rate limit of
$0.35$\,K per move. The predicted 20-bar pressure limit is softened with
penalty 25, and the objective combines terminal product with a gas-inventory
term of weight 0.003. An additional 1-bar opening margin is used for the
initial moves, except for PWL-native. For context, Table~\ref{tab:hydro-nb09} also retains
the full-remaining-horizon IPOPT runs, which instead use a horizon-averaged
inventory term of weight 0.05. Thus, those runs differ in objective treatment
as well as horizon. The reported solve times are the accumulated solver-call
times recorded in notebook 09; they exclude training, model construction,
and true-plant integration.

\begin{figure}
\centering
\includegraphics[width=\textwidth]{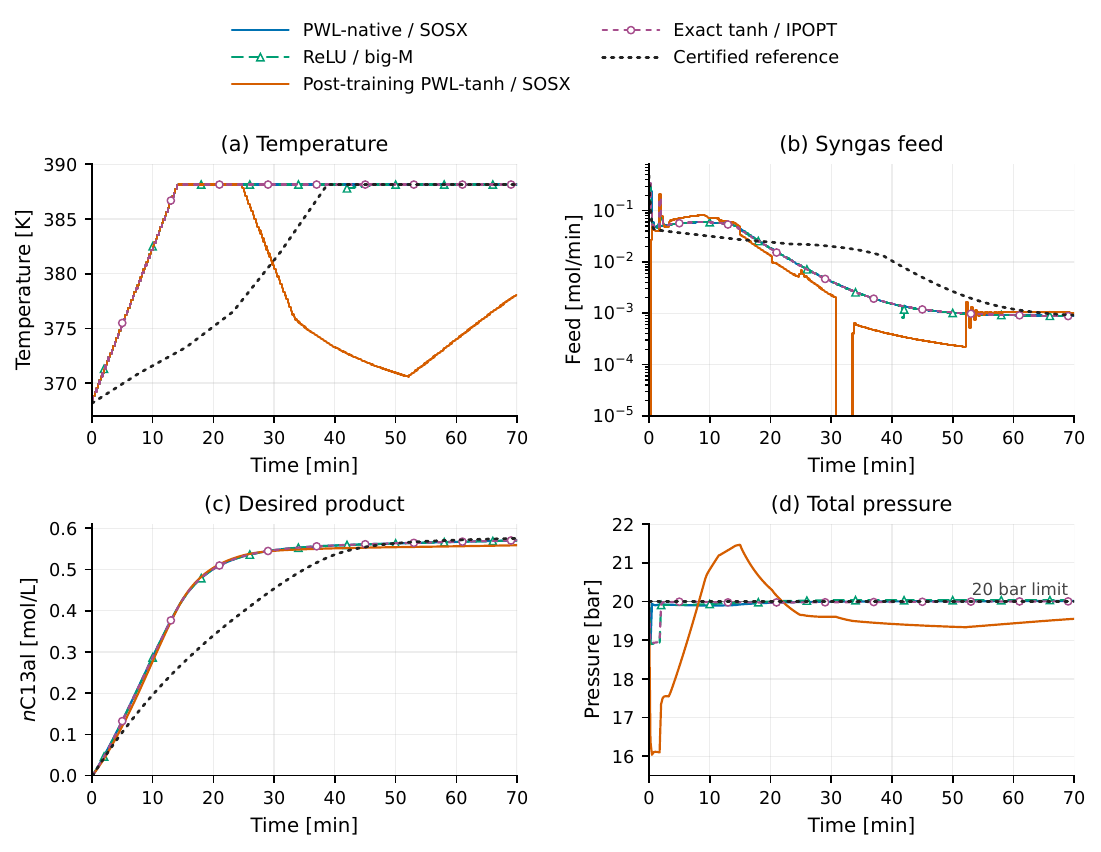}
\caption{Hydroformylation closed-loop comparison from notebook 09 over the
70-minute batch. The four controllers use $H=1$: PWL-native / SOSX,
ReLU / big-$M$, post-training PWL-tanh / SOSX, and exact tanh / IPOPT.
Temperature and syngas feed are applied to the true plant; product and
pressure are the resulting plant responses. The certified reference is
shown in black, with its 15-second control samples connected for visual
clarity. The feed axis is logarithmic, with zero feed displayed at
$10^{-5}$\,mol/min. The full pressure range of all four controllers is
shown; their precise peaks and solve times are given in
Table~\ref{tab:hydro-nb09}.}
\label{fig:hydro-pwl-native}
\end{figure}

Figure~\ref{fig:hydro-pwl-native} shows the four representative controllers,
and Table~\ref{tab:hydro-nb09} retains the plain, physics-informed, and pruned
variants. The one-step exact-tanh IPOPT controller reaches a terminal
$n$C13al of $0.5706$, approximately 99.2\% of the reference, with a maximum
realized pressure of $20.002$\,bar and a total solve time of $14.36$\,s.
It provides a favorable combination of product, pressure-limit adherence,
and runtime among the tanh-based controllers. ReLU / big-$M$ is faster at
$7.78$\,s, reaching $0.5703$ with a peak of $20.036$\,bar. PWL-native / SOSX
also reaches $0.5703$, with a peak of $20.013$\,bar, in $49.67$\,s. These
small excursions are reported explicitly rather than interpreted as strict
true-plant feasibility. Post-training PWL-tanh / SOSX gives a lower terminal
product of $0.5586$, a larger pressure peak of $21.466$\,bar, and a longer
solve time of $151.30$\,s. None of these short-horizon runs required a
fallback control move.

\begin{table}[htbp]
\caption{Notebook-09 closed-loop results for the 70-minute, 280-move batch
(seed 2479085). PWL-tanh denotes post-training activation replacement;
PI denotes physics-informed training. Full denotes the full remaining
horizon. The gap is reference minus realized terminal product, in mol/L;
product above the reference does not establish a feasible improvement.
Test MAE is the aggregate one-step state error of the smooth predictor for
tanh-trained networks and of the deployed predictor for PWL-native and
ReLU. Time is the recorded accumulated solve time for the batch.}
\label{tab:hydro-nb09}
\centering
\footnotesize
\setlength{\tabcolsep}{3pt}
\begin{tabular}{lccccccc}
\toprule
Controller & $H$ & Terminal & Gap & Max $p$ & Min $p$ & Test & Time \\
 & & $n$C13al & & [bar] & [bar] & MAE & [s] \\
\midrule
PWL-native / SOSX          & 1 & 0.5703 & 0.0051 & 20.013 & 18.896 & 0.00323 & 49.67 \\
PWL-tanh / SOSX            & 1 & 0.5586 & 0.0169 & 21.466 & 16.040 & 0.00363 & 151.30 \\
PI PWL-tanh / SOSX         & 1 & 0.5562 & 0.0192 & 20.000 & 16.834 & 0.13888 & 101.21 \\
Pruned PWL-tanh / SOSX     & 1 & 0.5646 & 0.0109 & 21.297 & 14.540 & 0.00577 & 91.96 \\
Exact tanh / IPOPT         & 1 & 0.5706 & 0.0049 & 20.002 & 18.934 & 0.00363 & 14.36 \\
PI exact tanh / IPOPT      & 1 & 0.5764 & $-0.0010$ & 20.601 & 19.109 & 0.13888 & 14.45 \\
Pruned exact tanh / IPOPT  & 1 & 0.5705 & 0.0049 & 20.091 & 18.677 & 0.00577 & 14.22 \\
Exact tanh / IPOPT         & Full & 0.5683 & 0.0072 & 21.663 & 16.382 & 0.00363 & 1273.26 \\
PI exact tanh / IPOPT      & Full & 0.5545 & 0.0210 & 20.576 & 17.048 & 0.13888 & 1206.55 \\
Pruned exact tanh / IPOPT  & Full & 0.5728 & 0.0026 & 20.076 & 19.658 & 0.00577 & 948.39 \\
ReLU / big-$M$            & 1 & 0.5703 & 0.0052 & 20.036 & 18.901 & 0.00263 & 7.78 \\
PI ReLU / big-$M$         & 1 & 0.5763 & $-0.0008$ & 20.617 & 18.966 & 0.13544 & 7.91 \\
Pruned ReLU / big-$M$     & 1 & 0.5701 & 0.0054 & 20.021 & 18.904 & 0.00372 & 8.05 \\
Certified reference      & -- & 0.5755 & 0 & 20.000 & 20.000 & -- & -- \\
\bottomrule
\end{tabular}
\end{table}

The full-horizon default-tanh IPOPT run reaches $0.5683$ with a pressure
peak of $21.663$\,bar in $1273.26$\,s, so its additional lookahead does not
improve on the one-step result for this network. Pruning improves that
full-horizon result to $0.5728$ and reduces the pressure peak to
$20.076$\,bar, but still requires $948.39$\,s. Physics-informed training
also has mixed effects: PI PWL-tanh / SOSX respects the pressure limit in
this run but reaches only $0.5562$, whereas PI ReLU and one-step PI tanh /
IPOPT exceed the reference product while also exceeding the pressure limit.
Neither a lower prediction error nor a larger terminal product alone is
therefore sufficient to rank these controllers.

A simple explanation is that the controller must predict the active
pressure constraint accurately where its own decisions take the plant.
Preserving the activation in PWL-native, ReLU, or exact-tanh embedding
removes activation-replacement error, but leaves the learned transition
error. A longer horizon can compound that error, and a soft predicted
pressure constraint does not enforce the true-plant limit. The gas feed
must consequently be judged by how well it maintains pressure near
20\,bar, alongside the terminal product. Differences in the temperature
profile are less conclusive: the reference's constraint-seeking feed arc
and sensitivity-seeking temperature arc can yield similar terminal
performance with different temperature histories
\citep{aydin2017dynamic,aydin2018nmpc}.

For this case, compact one-step controllers recover most of the reference
product at modest computational cost. Exact-tanh IPOPT gives a particularly
good balance of pressure adherence and runtime, ReLU / big-$M$ is the
fastest, and PWL-native / SOSX provides a consistent PWL representation
with comparable product recovery. The main lesson of the experiments is
to select the surrogate and horizon using realized plant constraints,
terminal performance, and solve time together; increasing the horizon,
adding physics penalties, or pruning is not uniformly beneficial.

\section{Conclusions}
\label{sec:conclusion}
This work presented PWLNet, an open-source Python framework for
constructing, compressing, reformulating, and validating neural-network
surrogates in nonlinear process optimization. The framework separates the
trained model from both the training framework and the optimization
backend through a framework-independent network representation, and
exposes ReLU big-M, exact smooth tanh, piecewise-linear tanh big-M/SOS2,
and SOSX branch-and-bound formulations, including MIQP and MISOCP
objective variants, behind one interface, together with
optimization-aware pruning, physics-informed training, interval bound
propagation, and process benchmark modules.

The four benchmark studies support several practical conclusions. First,
formulation choice is an empirical, application-specific decision: on
Williams--Otto the exact tanh NLP is the most accurate embedding but
offers no global certificate and may select slightly infeasible points,
while the PWL-tanh MILPs certify global optimality of the embedded model
in sub-second times at a modest profit cost dominated by the economics
relaxation and surrogate error rather than the breakpoint count.

PWL-native training provides consistency between the predictor fitted to
 data and the model embedded in optimization by retaining the same 
 piecewise-linear activation throughout the workflow. This eliminates 
 activation-replacement error, but does not eliminate prediction error 
 relative to the original process. The hydroformylation results demonstrate 
 the effectiveness of this approach, with the PWL-native controller achieving 
 a terminal product of $0.5703$, close to the reference value of $0.5755$, 
 while exhibiting only a minor pressure exceedance of $0.013$,bar. These results 
 demonstrate that piecewise-linear neural surrogates can recover near-optimal 
 solutions to the original nonlinear process while enabling globally certifiable 
 mixed-integer optimization, offering a potentially more tractable alternative to 
 nonlinear neural-network embeddings in applications where global optimization is challenging.

Second, optimization-aware compression works as intended: removing 15\% of the
hidden neurons of a ReLU surrogate cut the embedded MILP solve time by a
factor of 6.5 while leaving the recommended operating point essentially
unchanged, and the same qualitative size--tractability benefit held for
the recurrent surrogate. In the hydroformylation study, however, pruning
from 16 to 13 neurons increased the tanh test MAE from $0.0036$ to $0.0058$
and the ReLU test MAE from $0.0026$ to $0.0037$, showing that compression
is not uniformly accuracy-neutral. 

Third, constraint satisfaction must be
assessed in closed loop rather than inferred from aggregate surrogate error
alone. The Williams--Otto and hydroformylation results show that errors
in constrained outputs can drive true-plant violations. In hydroformylation,
physics penalties and pressure margins did not consistently improve the
realized policy, reinforcing the need for direct plant validation. 

Fourth, the
SOSX workflow provides globally optimal PWL-MAE training at small
incremental cost over a gradient warm start and a binary-free deployment
path that matches big-M solutions exactly, with incumbent warm starts and
persistent LP solves making it practical inside receding-horizon loops.

Finally, the dynamic case studies show that compact linearized surrogates
can approach reference terminal objective values in closed loop when the
prediction target and horizon handling encode the control task correctly,
as in the terminal-aware shrinking-horizon formulation required by the
inverse-response Van de Vusse problem. In Van de Vusse, the reported
99.2--99.7\% of the direct-shooting reference is achieved with state bounds
respected. In hydroformylation, one-step exact-tanh IPOPT achieves 99.2
the certified reference product in $14.36$,s, while ReLU achieves comparable 
performance in a shorter $7.78$,s. Although PWL-native requires longer solution times, 
it achieves nearly identical terminal performance to IPOPT, with an absolute objective 
gap of $0.0051$ compared with $0.0049$ for IPOPT, while enabling globally certifiable 
MILP optimization. All three approaches exhibit only minor pressure-limit exceedances, 
demonstrating that compact neural surrogates can recover near-optimal control solutions 
with small constraint deviations.

Future work includes adaptive breakpoint placement; reducing training complexity for single-layer MAE-PWL networks on large datasets; extending SOSX training beyond single-hidden-layer, fixed-readout networks; developing parametric and adaptive surrogates for robustness to plant drift; and establishing closed-loop guarantees for the embedded controllers.

\section{Software Availability}
\label{sec:software}

PWLNet is developed as an open-source Python package, with source code
available at \url{https://github.com/SunoySanyal/PWLNet}. The repository
is organized as a research codebase with a \texttt{pyproject.toml}
configuration, optional dependency extras for training and solver
backends, a test suite, and the fourteen executable benchmark
notebooks. The package is available on PyPI as \texttt{pwlnet}.

The repository is structured as follows:

\begin{verbatim}
PWLNet/
|-- README.md
|-- pyproject.toml
|-- requirements.txt
|-- docs/
|-- examples/
|-- notebooks/
|-- datasets/
|-- tests/
`-- src/
    `-- pwlnet/
        |-- core/
        |-- data/
        |-- training/
        |-- physics/
        |-- pruning/
        |-- formulations/
        |-- branching/
        |-- benchmarking/
        |-- objectives/
        |-- multiobjective/
        |-- solvers/
        |-- embedding/
        |-- recurrent/
        |-- processes/
        `-- utils/
\end{verbatim}

\section*{Data Availability}

The Boston housing dataset is publicly available through standard machine learning repositories. The Williams--Otto, Van de Vusse, and hydroformylation benchmark data is generated from the corresponding benchmark models and provided with the repository. All generated datasets, scripts, trained models, and optimization results are available in the project repository.

\section*{CRediT Author Statement}

Sunoy Sanyal: Conceptualization, Methodology, Software, Data curation, Validation, Formal analysis, Investigation, Visualization, Writing -- original draft.

Canberk Demir: Conceptualization, Methodology, Software, Data curation, Validation, Formal analysis, Investigation, Visualization, Writing -- original draft.

Metin T\"urkay: Supervision, Methodology, Resources, Project administration, Writing -- review and editing.

Erdal Ayd{\i}n: Supervision, Methodology, Resources, Project administration, Writing -- review and editing.

\section*{Declaration of Competing Interest}

The authors declare no known competing financial interests or personal relationships that could have appeared to influence the work reported in this paper. This research did not receive any specific grant from funding agencies in the public, commercial, or not-for-profit sectors.


\section*{Declaration of Generative AI and AI-Assisted Technologies in the Writing Process}

\paragraph{Writing Process}
During the preparation of this work the authors used ChatGPT (OpenAI) in order
to improve the language, grammar, and structure of author-written text.  The
tool was applied only to content that the authors had already produced
themselves: it corrected grammatical errors, improved sentence structure and
paragraph organisation, added a small number of connective sentences, and
removed redundant passages.  It was not used to generate scientific content,
data, results, figures, or references.

\paragraph{Software Development}
The accompanying \texttt{PWLNet} software was developed with the assistance of
Claude Code (Anthropic), an agentic coding assistant, which was used for
author-directed implementation, refactoring, testing, debugging, and
documentation of code written to author-specified designs and formulations.
All algorithms, model formulations, benchmark choices, and experimental
protocols were specified by the authors; all AI-suggested code was reviewed by
the authors before use; and every reported result was produced by executing the
authors' own code and verified against independent reference solutions.  The
complete source code is publicly available so that all results can be
reproduced and audited independently.

\paragraph{Responsibility}
After using these tools, the authors reviewed and edited the content as needed
and take full responsibility for the content of the publication.  No
generative-AI tool is listed as an author.

\section*{Funding}

This research did not receive any specific grant from funding agencies
in the public, commercial, or not-for-profit sectors.

\bibliographystyle{cas-model2-names}
\bibliography{references}

\end{document}